%% file: Main.tex
\documentclass[final, 12pt]{elsarticle}

\usepackage{graphicx}
\usepackage{amssymb}

\usepackage{lineno}

\include{Pack}

\journal{Applied Energy}

\pgfplotsset{compat=1.14}
\begin{document}

\begin{frontmatter}


\title{Flexible unit commitment of a network-constrained combined heat and power system}



\author[Skoltech]{Alvaro Gonzalez-Castellanos\corref{mycorrespondingauthor}}
\cortext[mycorrespondingauthor]{Corresponding author}
\ead{alvaro.gonzalez@skolkovotech.ru}

\author[UCDublin]{Priyanko~Guha Thakurta}
\author[Skoltech]{Aldo Bischi}

\address[Skoltech]{Center for Energy Systems - Skolkovo Institute of Science and Technology, Moscow, 143026, Russian Federation}
\address[UCDublin]{University College Dublin, Belfield, Dublin 4, Ireland}

\begin{abstract}
Large Combined Heat and Power (CHP) plants are often employed in order to feed district heating networks, in Europe, in post soviet countries and China. Traditionally they have been operated following the thermal load with the electric energy considered as a by-product, while the modern trend includes them in the electric market to take advantage of the flexibility they could provide. This implies the necessity to consider the impact on the electric grid while filling the thermal load requests.
%

A detailed Mixed Integer Linear Programming (MILP) optimization model for the solution of the network-constrained CHP unit commitment of the day-ahead operation is introduced. The developed model accounts for lossless DC network approximation of the electric power flow constraints, as well as a detailed characterization of the CHP units with useful effect, heat and power, function of one and two independent variables (\enquote{degrees-of-freedom}), and thermal energy storage. A computational validation of the outlined model on a CHP test system with multiple heating zones is presented in the form of computational test cases. The test cases illustrate the impact on the flexibility of the implementation of the energy storage, network constraints and joint multi-system operation. The conducted studies have highlighted the importance of a comprehensive and integrated analysis of multi-energy systems to exploit the operational flexibility provided by the cogeneration units. The  joint  operation  of  the  thermal  and  electric  system allows to reap economic, operational efficiency, and environmental benefits. The developed model can be easily extended to include diverse multi-energy systems and technologies, as well as more complex representations of the energy transmission networks, and the modeling of renewable energy resources dependent of one or more independent, weather-related, variables.\\
\end{abstract}


\end{frontmatter}

\begin{longtable}[!ht]{p{.175\textwidth}  p{.825\textwidth} }
\caption*{Nomenclature list}\\
\label{tab:Nomenclature}\\

	\hline
    \multicolumn{2}{l}{Sets} \\
    \hline
    $I$ &  Set of all generation units $i$ in the system\\
    $EP$ & Units that produce electric power\\
    $HT$ & Units that produce high-temperature heat\\
    $CHP$ & Units that produce both electricicity and heat ($EP \cap HT$)\\
    $I_1$ & Units with one independent variable\\
    $I_2$ & Units with two independent variables\\
    $J_i$ & Set of evaluation points to characterize the energy output based on the first independent variable for $i \in I_1$\\
    $K_i$ & Set of evaluation points to characterize the energy output based on the second independent variable for $i \in I_2$\\
    \\
    \multicolumn{2}{l}{Indexes}\\
    \hline
    $i$ & Generation unit\\
    $j$ & Sampling point for the first independent variable\\
    $k$ & Sampling point for the second independent variable\\
    $n$ & Power node. $w = alias(n)$\\
    $t$ & Time step. $\tilde{t}= alias(t)$\\
    $z$ & Thermal zone\\
    \\ 
    \multicolumn{2}{l}{Superscripts}\\
    \hline
    d & Demand\\
    f & Fuel\\
    low & Related to the lower triangle for bi-dimensional piecewise linearization\\
    OM & Operation \& Maintenance\\
    s & Storage-related\\
    st & Startup\\
    tm & Time-related\\
    up & Related to the upper triangle for bi-dimensional piecewise linearization\\
    \\
    \multicolumn{2}{l}{Parameters}\\
    \hline
    $B_{nw}$ & Susceptance between the nodes $w$ and $n$, p.u.\\
    $C_i^\textup{f}$ & Cost of fuel for unit $i$, \euro \\
    ${C}^\textup{OM,f}_i$ & Operation \& Maintenance cost by fuel consumption of the unit $i$, \euro \\
    ${C}^\textup{OM,st}_i$ & The cost of Operation \& Maintenance of the unit $i$ linked to its start up procedures, \euro\\
    $C^\textup{OM,tm}_i$ & The cost of Operation \& Maintenance of the unit $i$ linked to its operational time, \euro\\
    $C^\textup{st}_i$ & Cost penalization per startup procedure for $i$, \euro\\
    $DR_i$ & Down-ramp rate for unit ${i \in EP}$, MWh\\
    $\hat{E}_{i,j}$ & Electric energy generation of unit $i \in EP \cap I_1$, correspondent to the evaluation point $j$, MWh\\
    $\hat{E}_{i,j,k}$ & Electric energy generation of ${i \in EP \cap I_2}$, correspondent to the evaluation point ${(j,k)}$, MWh\\
    $\overline{E}_i$ & Maximum electricity output for ${i \in EP}$, MWh\\
    $\underline{E}_i$ & Minimum electricity output for ${i \in EP}$, MWh\\
    $E^\textup{d}_{n,t}$ & Electric energy load at $n$ during time $t$, MWh\\
    $\hat{F}_{i,j}$ & Evaluation point $j$ for the fuel consumption (first independent variable) of the unit $i$, MWh\\
    $\overline{\Phi}_{kw}$ & Maximum permissible active power through power line connected between $n$ and $w$, MW\\
    $\Gamma_{i,z}$ & Binary parameter representing if the unit $i \in HT$ serves the zone $z$\\
    $\eta_z$ & Efficiency for the transfer of heat within zone $z$, \%\\
    $\hat{H}_{i,j}$ & Thermal energy generation of $i \in HT \cap I_1$, correspondent to the evaluation point $j$, MWh\\
    $\hat{H}_{i,j,k}$ & Thermal energy generation of ${i \in HT \cap I_2}$, correspondent to the evaluation point $(j,k)$, MWh\\
    $\overline{H}_i$ & Maximum heat output of unit $i \in HT$, MWh\\
    $\underline{H}_i$ & Minimum heat output of unit $i \in HT$, MWh\\
    $H^\textup{d}_{z,t}$ & Total thermal load for $z$ during time $t$, MWh\\
    $L_z^{s}$ & Percentage of the storage level that is lost per hour in the storage of the zone $z$, \% \\
    $N_i^\textup{st}$ & Maximum number of start up processes tolerated by the unit $i$\\
    $MT_i$ & Minimum time that the unit $i \in I$ must be kept turned on, $h$\\
    $\hat{O}_{i,k}$ & Break-point $k$ for the second independent variable of the unit $i \in I_2$\\
    $\Omega_{in}$ & Binary parameter representing if the electric unit $i \in EP$ is connected to node $n$\\
    $\hat{P}_{i,j}$ & Electric power generation of unit $i \in EP \cap I_1$, correspondent to the evaluation point $j$, MW\\
    $\hat{P}_{i,j,k}$ & Electric power generation of unit $i \in EP \cap I_2$, correspondent to the evaluation point ($j,k$), MW\\
    $P^\textup{d}_{n,t}$ & Electric power load at $n$ during $t$, MW\\
    $\overline{S}_z$ & High-temperature heat storage capacity of the zonal storage for zone $z$, MWh\\
    $UR_i$ & Up-ramp rate for unit ${i \in EP}$, MWh\\
    
    \\
    \multicolumn{2}{l}{Variables}\\
	\hline
    $\alpha_{i,j,t}$ & Variable associated with each breakpoint $j$ of the unit $i \in I_1$ during $t$.\\
    $\alpha_{i,j,k,t}$ & Variable associated with each breakpoint $(j,k)$ of the unit $i \in I_2$ during $t$.\\
    $\delta_{n,t}$ & Phase angle of node $n$ during $t$.\\
    $e_{i,t}$ & Electric energy generation of $i \in EP$ during $t$.\\
    $f_{i,t}$ & Fuel consumption of unit $i$ during $t$.\\
    $h_{i,t}$ & Total heat production of unit $i \in HT$ during $t$.\\
    $h_{i,z,t}$ & Thermal energy generation of unit $i \in HT$ in zone $z$ during $t$.\\
    $p_{i,t}$ & Electricity production of unit $i \in EP$ during $t$.\\
    $p_{nw,t}$ & Active power flowing through power line connected between nodes $n$ and $w$ during $t$.\\
    $s_{z,t}$ & Level of the thermal storage in zone $z$ during $t$.\\
    $\beta_{i,j,t}$ & Dummy binary variable associated with the characterization of the performance curves of the unit $i \in I_1$ in the $j^\textup{th}$ interval $[f_j,f_{j+1}]$ during $t$.\\
    $\beta_{i,j,k,t}^\textup{up}$ & Binary variable associated with the upper triangle of the rectangle corresponding to the intervals $[f_j,f_{j+1}]$ and $[o_k,o_{k+1}]$, for the characterization of the performance curves of $i \in I_2$ during $t$.\\
    $\beta_{i,j,k,t}^\textup{low}$ & Binary variable associated with the lower triangle of the rectangle corresponding to the intervals $[f_j,f_{j+1}]$ and $[o_k,o_{k+1}]$, for the characterization of the performance curves of $i \in I_2$ during $t$.\\
    $\tau_{i,t}$ & Binary variable symbolizing if $i$ has been turned on during $t$.\\
    $\theta_{i,t}$ & Binary variable representing if $i$ is turned on or off during $t$.\\
    
    \hline

\end{longtable}

\section{Introduction} \label{S: Intro}
Large-scale Combined Heat and Power (CHP) plants, above 10MW, are being increasingly deployed within the industrial and public sectors. By 2011, their share increased to 79~\% of the total thermal energy in District Heating (DH) networks, thereby making it an essential element in the energy mix of countries that endure long winters \cite{FernandezPales2014}. The reason behind such an increase is two-fold: due to their higher efficiency (first principle thermodynamics efficiency) in converting the primary energy, fuel, into heat and power (useful effect) compared to separate conventional generation, and their flexibility in varying the share of generated heat and power depending on the adopted technology, as is the case of the extraction condensing steam turbine. The importance of the heat produced by the CHP plants has been one of the leading arguments for the introduction of DH \cite{Werner2017a}. This has also been reflected on their share of generated electric power, reaching as much as 50~\% of the total electricity generation in countries such as China, Latvia, Russia, Finland, and Denmark \cite{Kerr2009}.

CHP units have their heat and electric power generation depending on the number of independent variables - i.e., \enquote{degrees-of-freedom}. The so-called \enquote{one-degree-of-freedom} units consider the load percentage representing the consumed fuel as the one and only independent variable e.g. simple cycle gas turbine, whereas other units can be controlled based on two independent variables, \enquote{two-degrees-of-freedom}, which allow to decouple the heat and electricity generation e.g. extraction condensing steam turbine whose independent variables are the above-mentioned fuel and the valve opening controlling the ratio of heat and electric power. Such generation flexibility could allow for the balancing of fluctuating renewable energy resources. The CHP units could provide a higher thermal output when the renewable generation is higher, and contribute to the electric energy balance once the renewable resources contribution diminishes.

The current practice of operating CHPs is a decoupled one. In other words, the units are primarily following the thermal load. After the thermal working point has been set on the comprehensive generation curve, an electric economic dispatch is performed. The generated electricity is thus constrained by the generated heat, and its generation model is simplified accordingly. Such an approach does not achieve a global minimum cost of system operation, since the optimization of the thermal and electric systems are performed in a sequencial, rather than integrated way. While accounting for energy storage, a decoupled dispatch diverges further from achieving the minimum operating cost. This corresponds to the fact that when combined with CHP units, the use of energy storage facilitates load shifting, increasing the flexibility of the system \cite{Lund2015, NordicTSOs2017}.\\
Moreover, the electric power flow constraints must also be accounted for. The lack of consideration of the electric power flow could lead to safety violations, resulting in the loss of served load. The average annual value of the economic losses due to unserved load is 10,000~\euro/MWh and its value is higher during the winter months for households of EU Member States in Northern Europe \footnote{This value is estimated based on the amount that the end users would be willing to pay to their energy retailer to guarantee an uninterrupted electricity supply \cite{Vassilopoulos2003}.}.Hence, a coupled heat and electric power dispatch of CHPs, while accounting for electric transmission flow constraints, becomes economically attractive \cite{Shivakumar2017a}.

\subsection{Relevant literature}
The non-linear non-convex nature of energy generation curves of CHP units along with the electric power flow equations makes the coupled CHP unit commitment and economic dispatch (CHPED) problem non-linear and non-convex. Moreover, the binary states associated with on/off and start up of the units result in a Mixed Integer Non-linear Programming (MINLP) CHPED model. Solving such a computationally challenging problem \cite{Belotti2013} is handled by piecewise linearization of the performance curves which, in turn, transform the problem into a Mixed Integer Linear Programming (MILP) one \cite{Bischi2014c}. This formulation takes advantage of the effectiveness of state-of-the-art linear solvers instead of nonlinear ones \cite{Taccari2015}. Piecewise linearization approach has proven to be close to the real solution with a moderate amount of piecewise intervals and has already been implemented in solving CHP dispatch problems \cite{Bischi2014c}. Moreover, mathematical \cite{Rooijers1994a, Sashirekha2013a} and heuristic methods \cite{Mellal2015a,Basu2015a,Jayakumar2016a} are also employed in solving such a problem.

Several literature have addressed the integration of electric power flow in integrated heat and power models. Among them are the methods based on the development of non-convex models for the electric and thermal networks \cite{Liu2016c}, and the  integration of gas networks \cite{Liu2016}.\\
Additional works have focused on the development of convex models for CHP systems based on multi-energy virtual power plants \cite{Capuder2016a}, i.e., distributed energy generation operated as one larger plant; and integration of renewable energy resources: wind \cite{Saint-Pierre2016a,Lu2013}, and solar thermal \cite{Wang2015}. The focus of these convex models was mainly on energy flows between the systems neglecting the effect of the network constraints in the flexible operation of the units. In the previously described works the electric, thermal and CHP units are modeled as constant efficiency units, thereby leading to a strongly simplified characterization given the nonlinear nature of the performance curves of the units. An accurate representation of the energy generation at partial load becomes essential when considering short-term scenarios and assessing the system flexibility \cite{Abdin2018}. This corresponds to the fact that in a constant-efficiency-based model the generation units do not posses a measurement of how adjustable the power and heat production of a CHP can be, allowing it to produce more energy for the daily electricity demand peaks when the thermal request is low, and viceversa. Conducing to an over sizing of the required generation capacity.\\
Rong \textit{et al.} \cite{Rong2016d} developed a model for the optimization of a CHP dispatch with multiple generation and consumption sites based on a dynamic programming model. In this study, the electric power flow is calculated based on an energy flow model, without the inclusion of Kirchhoff’s voltage law. Hence, the electric power was transmitted between the nodes based on a nodal power balance, surplus and shortage, rather than based on the actual path that the current would follow based on the electric characteristics of the transmission line. This could lead to the overestimation of the amount of power flowing between two nodes and its related generation output, making it necessary to redispatch some of the scheduled generators and incurring in scheduling compensation costs. The balance for the modeled heat sites was based on a load-generation balance, the use of thermal energy storage was not included.

An optimization framework for the integration of thermal and electric energy system while assessing the network constraints of the electric grid is defined by Morvaj \textit{et al.} \cite{Morvaj2016a}. The power flow is successfully integrated within the CHPED problem in a radial distribution network. 

Nuytten \textit{et al.} \cite{Nuytten2013} evaluated the increase in operation flexibility consequence of the combination of a CHP unit with thermal energy storage. In this work, the maximum flexibility of the system as a function of the CHP and storage capacity is evaluated. The modeled CHP plant had its electricity and heat generation as a function of the consumed fuel, i.e. \enquote{one-degree-of-freedom}. A linear relationship between the size of the energy storage and the available flexibility was observed, whereas the flexibility saturates with the increase in CHP capacity, i.e., once the rated thermal capacity of the CHP plants reaches a fourth of the peak thermal demand, further increase in installed capacity will not improve the system flexibility. The flexibility gains for the coupling of the CHP unit with an electric system with photo-voltaic generation (PV) are also analyzed. The inclusion of the CHP with energy storage allowed the system to reduce its electricity export around noon and import during the night. This was achieved by delaying the CHP energy generation to the evenings, when the PV generation diminished. The mismatched heat demand was fulfilled, during these hours, by the thermal energy storage. For the flexibility analysis only an energy balance is considered. Thus, being necessary an assessment of the flexibility gains derived from the employment of CHP units alongside thermal energy storage in a network constrained multi-energy system.

\subsection{Paper Contributions and Organization}
The main contribution of this paper is the formulation of a model that allows the assessment of the flexibility gains derived from the integrated operation of a combined heat and power system. For this purpose, it is devised a combined heat and power unit commitment that allows to introduce electric network constraints in the operation of CHP units with multiple degrees-of-freedom. The developed model includes thermal energy storage and district heating, along with start up costs and ramp constraints. A detailed characterization of the generation units is implemented through the piecewise linearization of their useful effect as a function of one and two variables (degrees-of-freedom). The implemented unit modeling and a network-constrained electric system allows for an analysis of the impact of the CHP flexibility on the operation of an electric transmission system. The proposed model is a Mixed Integer Linear Program (MILP) one.

Given the higher prices for electric energy storage compared to its thermal counterpart \cite{Lund2016a}, the developed framework exploits the thermal storage in order to store the heat from CHPs when the electric power transmission system becomes constrained. The sizing of the thermal energy storage follows the practices employed in modern district heating systems, such as the one in the city of Turin, where the district heating network has a storage system with a total volume of more than 12,000~cubic meters\footnote{The system stores around 10~\% of the annual energy consumed in the city, which is about 178~GWh of thermal energy.} \cite{Baccino2014a}.\\
Several configurations of CHP systems are considered with the use of 4 test cases, which focus on: unit flexibility; storage integration; network constraints; and coordination of the electric and thermal system.\\
The paper is organized as follows: Section \ref{S: MathModel} proposes the mathematical formulation to include CHPs within the unit commitment and economic dispatch optimization framework. The results of the proposed framework are shown in Section \ref{S: Test}. Finally, Section \ref{S: Conclusion} draws the conclusions of the paper.

\section{Mathematical Model} \label{S: MathModel}
\subsection{Objective function}
For a daily operation, only the operational conditions of the system are considered. Long-term expenditures such as installation, amortization and legal costs, are not introduced in the analysis. Therefore, the objective of the optimization model is to determine the minimum operational cost of the joint energy system while satisfying the thermal and electrical constraints, and is formulated as:
\begin{IEEEeqnarray}{rl} \label{Eq: Obj}
    \text{min}\quad & \sum_{t}\Big[{C_{t}^\textup{f}} + {C_{t}^\textup{st}} + {C_{t}^\textup{OM}}\Big] \IEEEnonumber\\
    \text{subject to:}\quad & \text{Scheduling equations \eqref{Eq: Total Fuel} to \eqref{Eq: Lim S}}
\end{IEEEeqnarray}

The operational costs of the system can be divided in:
\begin{itemize}
    \item Fuel consumption: 
        \begin{IEEEeqnarray}{C} \IEEEyessubnumber \label{Eq: Total Fuel}
            C_{t}^\textup{f} = \sum_i \hat{c}_{i}^\textup{f} \cdot f_{i,t}, \qquad \forall t.
        \end{IEEEeqnarray}
    \item Start-costs: during the start up procedure, due to the thermal inertia of the plant components, the plant consumes primary energy without producing useful effect. Thus, incurring in an extra cost of fuel consumption, which is accounted by
        \begin{IEEEeqnarray}{C} \IEEEyessubnumber \label{Eq: Total St}
            C_{t}^\textup{st} = \sum_i \hat{c}_{f,st,i} \cdot \tau_{i,t}, \qquad \forall t.
        \end{IEEEeqnarray}
    \item Operation and Maintenance (O\&M): the O\&M cost of each unit depends on the amount of time that the unit runs, the number of start up procedures that it undergoes and its primary energy consumption. It is given by
        \begin{IEEEeqnarray}{C} \IEEEyessubnumber \label{Eq: Total OM}
            C_{t}^\textup{OM} = \sum_i (\hat{c}_{i}^\textup{OM,tm} \cdot \theta_{i,t} +  \hat{c}_{i}^\textup{OM,st} \cdot \tau_{i,t} + \hat{c}_{i}^\textup{OM,f} \cdot f_{i,t} ), \qquad \forall t.
        \end{IEEEeqnarray}
\end{itemize}

\subsection{Energy systems modeling}
\subsubsection{Electric energy system}
The electric energy system can be modeled by
\begin{IEEEeqnarray}{rCll}
\sum_{n} \Big [ E^{\text{d}}_{n,t} - \sum_{i \in EP} \Omega_{in} \cdot e_{i,t} \Big ] &=& 0, & \qquad  \forall t \IEEEyesnumber \IEEEyessubnumber* \label{Eq: E Balance}\\
\sum_{n} \Big [ P_{n,t}  - \sum_{i \in EP} \Omega_{in} \cdot p_{i,t} - \sum_{w} p_{nw,t} \Big] &=& 0, & \qquad \forall t \label{Eq: Power Flow}\\
B_{nw}\cdot(\delta_{n,t}-\delta_{w,t}) &=& p_{nw,t}, & \qquad \forall t \label{Eq: LineFlow}\\
|p_{nw,t}| & \leq & \overline{\Phi}_{nw}, & \qquad \forall t.	\label{Eq: STR}
\end{IEEEeqnarray}

The electric energy balance is given by \eqref{Eq: E Balance}. The power balance at a node based on a DC power flow is given by \eqref{Eq: Power Flow}. The Static Thermal Rating (STR) is defined as the maximum permissible current through a line. \eqref{Eq: LineFlow} represents the power flow through a transmission line, whereas its STR is given by \eqref{Eq: STR}.

\subsubsection{Thermal energy system}
The thermal energy system is modeleded by
\begin{IEEEeqnarray}{rClr}
    H^\text{d}_{z,t} &\geq& \eta_z \Bigg[ \sum_{i \in HT} \Gamma_{iz}  h_{i,z,t} + (s_{z,t} - s_{z,t+1}) \Bigg] -{L_z^\text{s} \cdot{s_{z,t}}} ,&\quad \forall z,t\hspace{0.05cm} \IEEEyesnumber \IEEEyessubnumber* \label{Eq: H Balance}\\
    h_{i,t} &=& \sum_{z} \Gamma_{i,z} h_{i,z,t},& \quad \forall i \in HT, t. \label{Eq: Sum h}
\end{IEEEeqnarray}

The total thermal load is considered per zone, which is represented by the variable $\hat{H}_{Load,z,t}$, while the level of the zonal thermal storage at time $t$ is given by $s_{HT,z,t}$. The thermal energy balance for each zone $z$ of the system is given by \eqref{Eq: H Balance}. The inequality in \eqref{Eq: H Balance} corresponds to the possibility of the CHP units to dissipate excess heat into the environment, if it is economically advantageous. Since the thermal system is divided into geographical zones, it must be noted that the manipulation of the binary parameter $\Gamma_{i,z}$ would allow the unit $i \in HT$ to serve one or more thermal zones, depending on the topology of the system.The sum of the heat transferred into one or multiple regions cannot exceed the technical limits of the unit. The total thermal energy produced by the unit $i \in HT$ in the time t is given by \eqref{Eq: Sum h}.

The losses in the thermal system can be grouped as follows:
\begin{itemize}
\item Losses in the heat distribution network: they account for 8-10~\% of the total transferred heat \cite{Comakli2004a}. The constant $\eta_{z}$ represents the efficiency of the heat transfer within the zone $z$. A conservative value of 92~\% is assumed for $\eta_{z}$. 
\item Losses in the heat storage system: they can be considered as a fixed percentage of the storage level for each hour, accounted by the parameter $L_z^\text{s}$. This parameter depends on the technical characteristics of the storing device. The efficiency assumed for the storage units is of 98~\% \cite{Bischi2014c}, making the parameter $L_z^\text{s}$ equal to 2~\%.
\end{itemize}

\subsection{Unit characterization}
\subsubsection{Start up procedures}
Equations \eqref{Eq: On1}-\eqref{Eq: On3} define the value of the variable $\tau_{i,t}$, which represents the undergoing of a start up procedure, thus ensuring that it equals to 1 only for the time step when the unit $i$ is turned on.
\begin{IEEEeqnarray}{rCCCll}
    \tau_{i,t} &\leq& \theta_{i,t}, &&&\qquad\forall i, t	\IEEEyesnumber \IEEEyessubnumber* \label{Eq: On1}\\
    \tau_{i,t} &\leq& 1 &-& \theta_{i,t-1}, &\qquad\forall i, t	\label{Eq: On2}\\
    \tau_{i,t} &\geq& \theta_{i,t} &-& \theta_{i,t-1}, &\qquad \forall i, t.	\label{Eq: On3}	
\end{IEEEeqnarray}

The lifetime of a plant is effectively reduced by the number of start up and shutdown procedures that it undergoes. This is a consequence of the high levels of mechanical stress imposed in the prime mover by this dynamic behavior. Therefore, in order to maximize its lifetime, a plant cannot be submitted through more than an established number of start ups in a given period. Constraint \eqref{Eq: Max St} sets the maximum number of start up procedures per unit during a day
\begin{IEEEeqnarray}{C}
    \sum_{t} \tau_{i,t} \leq N_{i}^\textup{st}, \qquad \forall i.	\IEEEyessubnumber \label{Eq: Max St}
\end{IEEEeqnarray}

The minimum period that a unit remains committed is given by
\begin{IEEEeqnarray}{C}
    \theta_{i,t} \geq \sum_{\tilde{t}=1}^{MT_i} \tau_{i,t-\tilde{t}}, \qquad \forall i, t.	\IEEEyessubnumber* \label{Eq: Min Time}
\end{IEEEeqnarray}

\subsubsection{Performance curves for units with one independent variable}

A nonlinear dependency is present between the primary energy consumption and the useful effect production, heat, and electricity, in the generation units. Therefore, it is necessary to model the characteristic curves as functions that guarantee their convexity and that of the system. For this purpose, the performance curves of the units are characterized using a piecewise linear approximation, using the one-dimensional method presented by D’Ambrosio et al. \cite{DAmbrosio2010a}. To characterize the units, their characteristic curves are sampled through $J_i$ breakpoints. Expressions \eqref{Eq: Sum Alpha-Beta 1}-\eqref{Eq: H_1} allow the modeling of the performance curves as a piecewise linear function.
\begin{IEEEeqnarray}{rCll}
    \alpha_{i,j,t} &\leq& \beta_{i,j-1,t} + \beta_{i,j,t}, &\qquad \forall i \in I_1, j, t	\IEEEyesnumber \IEEEyessubnumber* \label{Eq: Sum Alpha-Beta 1}\\
    1 &=& \sum_{j=1}^{|J_i|-1} \beta_{i,j,t}, &\qquad \forall i \in I_1, t	\label{Eq: Sum Beta 1}\\
    1 &=& \sum_{j} \alpha_{i,j,t}, &\qquad \forall i \in I_1, t \label{Eq: Sum Alpha 1}\\
    f_{i,t} &=& \sum_{j}\alpha_{i,j,t} \cdot  \hat{F}_{i,j}, &\qquad \forall i \in I_1, t	\label{Eq: F_1}\\
    e_{i,t} &=& \sum_{j}\alpha_{i,j,t} \cdot  \hat{E}_{i,j}, &\qquad \forall i \in EP \cap I_1, t	\label{Eq: E_1} \\
    p_{i,t} &=& \sum_{j}\alpha_{i,j,t} \cdot  \hat{P}_{i,j}, &\qquad \forall i \in EP \cap I_1, t	\label{Eq: P_1}	\\
    h_{i,t} &=& \sum_{j}\alpha_{i,j,t} \cdot  \hat{H}_{i,j}, &\qquad \forall i \in {HT \cap I_1}, t.	\label{Eq: H_1}
\end{IEEEeqnarray}

Expressions \eqref{Eq: Sum Alpha-Beta 1}-\eqref{Eq: Sum Alpha 1} sets the approximation of the characterized function to be a convex combination of the extremes of the interval of interest. The fuel consumption of a unit with one degree of freedom is given by \eqref{Eq: F_1}. Constraints \eqref{Eq: Sum Alpha-Beta 1}-\eqref{Eq: P_1} apply to the units $i \in EP$ that produce electric power; while \eqref{Eq: Sum Alpha-Beta 1}-\eqref{Eq: F_1}, \eqref{Eq: H_1} model the units $i \in HT$ producing high temperature heat. CHP units are modeled by employing \eqref{Eq: Sum Alpha-Beta 1}-\eqref{Eq: H_1}.\\

\subsubsection{Performance curves of units with two independent variables}
The units with two degrees of freedom are those whose useful effect, heat and electricity, depend on the value of two independent variables, $f_{i,t}$ and $o_{i,t}$. Their characteristic curves are sampled through $J_i$ and $K_i$ break-points for the independent variables $f_{i,t}$ and $o_{i,t}$, respectively. The modeling of the characteristic curves of the units with two degrees of freedom is done following the triangle method described by D'Ambrosio et al \cite{DAmbrosio2010a}. The constraints \eqref{Eq: Sum Alpha-Beta 2}-\eqref{Eq: H_2} allow the modeling of the performance curves of the different types of units as a piecewise linear function.\\

\begin{figure}
    \centering
    \includegraphics{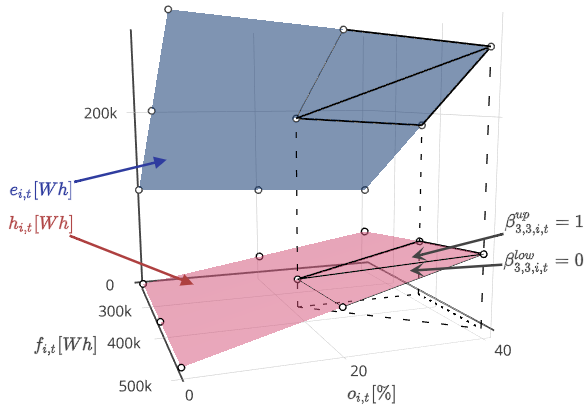}
    \caption{Triangle method for a CHP unit with two degrees-of-freedom}
    \label{Fig: PWL}
\end{figure}

\begin{IEEEeqnarray}{rClr}
    \alpha_{i,j,k,t} &\leq& \beta_{i,j,k,t}^\textup{up} + \beta_{i,j,k-1,t}^\textup{up} + \beta_{i,j-1,k-1,t}^\textup{low}  \IEEEnonumber \\
    &+& \beta_{i,j,k,t}^\textup{low}  + \beta_{i,j-1,k,t}^\textup{low} + \beta_{i,j-1,k-1,t}^\textup{up}, & \quad \forall i \in I_2, t \IEEEyesnumber \IEEEyessubnumber* \label{Eq: Sum Alpha-Beta 2}\\
    1 &=& \sum_{j}^{|J_i|-1} \sum_{k}^{|K_i|-1} (\beta_{i,j,k,t}^\textup{up} + \beta_{i,j,k,t}^\textup{low}), &\quad \forall i \in I_2, t \label{Eq: Sum Beta 2}\\	
    1 &=& \sum_{j} \sum_{k} \alpha_{i,j,k,t}, &\quad \forall i \in I_2, t \label{Eq: Sum Alpha 2}\\
    f_{i,t} &=& \sum_{j} \sum_{k}  \alpha_{i,j,k,t} \cdot \hat{F}_{i,j}, &\qquad \forall i \in I_2, t \label{Eq: F_2}\\
    o_{i,t} &=& \sum_{j} \sum_{k}  \alpha_{i,j,k,t} \cdot \hat{O}_{i,k}, &\qquad \forall i \in I_2, t \label{Eq: O_2}\\
    e_{i,t} &=& \sum_{j} \sum_{k}  \alpha_{i,j,k,t} \cdot \hat{E}_{i,j,k}, &\qquad \forall i \in EP \cup I_2, t \label{Eq: E_2}\\
    p_{i,t} &=& \sum_{j} \sum_{k}  \alpha_{i,j,k,t} \cdot \hat{P}_{i,j,k}, &\qquad \forall i \in EP \cup I_2, t \label{Eq: P_2}\\
    h_{i,t} &=& \sum_{j} \sum_{k}  \alpha_{i,j,k,t} \cdot \hat{H}_{i,j,k}, &\qquad \forall i \in HT \cup I_2, t. \label{Eq: H_2}
\end{IEEEeqnarray}
Figure \ref{Fig: PWL} displays the implementation of the triangle method on a CHP unit with a back-pressure steam turbine \cite{Bischi2014c}. The unit is characterized through 3 sampling points for the fuel - $f_{i,t}$ and the valve opening - $o_{i,t}$. As seen in the figure, a greater fuel input increases both electric - $p_{i,t}(f,o)$ and thermal generation - $h_{i,t}(f,o)$, respectively represented by the blue and red areas. Whereas, a bigger opening of the valve increases the thermal generation, while reducing the electricity output. In this example the CHP is operating on the upper triangle of the region delimited by the second and third sampling point.

\subsubsection{Technical limits for generation and storage units}
The economic parameters and the technical limitations provide an operational range for the unit $i$, this limits are represented by \eqref{Eq: Lim E} and \eqref{Eq: Lim H}
For the thermal storage, the amount of heat accumulated must be non-negative and below the maximum technical limit of the storage unit. These conditions are introduced for the zonal storage $z \in ZTh$ by constraint \eqref{Eq: Lim S}.\\
\begin{IEEEeqnarray}{rCCClr}
    \theta_{i,t}\underline{E}_{i} &\leq& e_{i,t} &\leq& \theta_{i,t}\overline{E}_{i}, & \forall i \in EP, t \IEEEyesnumber \IEEEyessubnumber* \label{Eq: Lim E}\\
    \theta_{i,t}\underline{H}_{i} &\leq& h_{i,t} &\leq& \theta_{i,t}\overline{H}_{i}, & \qquad \forall i \in HT, t\label{Eq: Lim H}\\
    0 &\leq& s_{z,t} &\leq& \overline{S}_{z}, & \forall z, t. \label{Eq: Lim S}
\end{IEEEeqnarray}
Ramp rates limit the power increase and decrease in consecutive periods of time, here represented by an ramp-up rate, $UR_i$, and ramp-down rate, $DR_i$, respectively \cite{Jayakumar2016a}. The ramp rates are modeled by
\begin{IEEEeqnarray}{rCCClr}
    -DR_i &\leq& e_{i,t} - e_{i,t-1} &\leq& UR_i, & \qquad \forall i \in EP, t.  \label{Eq: Ramp}
\end{IEEEeqnarray}

\section{Test system} \label{S: Test}
In order to demonstrate the effectiveness of the developed optimization model, a test case made of an electric transmission network and two heating zones is devised and then analyzed under four operational conditions. The day-ahead scheduling of the CHPs is done for 24 hours at an hourly time step. The hourly load at each node of the electric network, as well as the cumulative load at each thermal zone are given. Techno-economic parameters of the system and units are based on \cite{Bischi2014c}, and given in Table \ref{tab1}.\\
\setcounter{table}{0} 
\renewcommand{\arraystretch}{1.25}
\begin{table}[htb]
  \centering
  \caption{Techno-economic parameters}
    \begin{tabular}{p{5.5em}p{15em}c}
    \toprule
    \toprule
    Parameter & Basis & \multicolumn{1}{p{5em}}{Value} \\
    \midrule
    $C_i^\textup{f}$ & Thermal energy, LHV basis [\euro/kWh] & 0.6 \\
    $C^\textup{st}_{i}$ & Maximum energy input [\euro/kWh] & 0.009 $\cdot \hat{F}_{i,|J_i|}$ \\
    $C^\textup{OM,f}_i$ & Energy input [\euro/kWh] & 0.001 \\
    $C^\textup{OM,tm}_i$ & Operating hours [\euro/h] & 0.001 \\
    $C^\textup{OM,st}_i$ & Number of startups [\euro/on] & 1 \\
    $\eta_z$     & Thermal distribution efficiency [\%] & 92 \\
    $L^\textup{s}_z$  & Storage losses [\%] & 2 \\
    $|J_i|$,  $|K_i|$ & Sampling intervals for piecewise linearization of the units & 3 \\
    $N^\text{st}_i$ & Maximum number of startup procedures for a 24h period & 2 \\
    $MT_i$ & Minimum online time & 2 \\
    \bottomrule
    \bottomrule
    \end{tabular}%
  \label{tab1}%
\end{table}%

\subsection{Electric system}
The electric system is a modified IEEE 30-Bus Test Case \cite{Christie1993}. Two CHP plants CHP-1d and CHP-2d have been added to the standard model, located at nodes 5 and 2, and thermal zones 2 and 1 respectively. The CHP-1d is a gas turbine with heat recovery while the CHP-2d is a natural gas combined cycle with back-pressure steam turbine, i.e., its two independent variables are fuel consumption - $f_{i,t}$, and valve opening - $o_{i,t}$. Table \ref{Tab: EP Data} presents the technical parameters for the electric generation units. Their connection bus, up ramp limit - $UR$, down ramp limit - $DR$, as well as their sampling points for fuel consumption - $\hat{F}_{i,j}$ and electricity generation - $\hat{P}_{i,j}$ are presented. Table \ref{Tab: CHP1d Data} provides the technical parameters for the CHP unit with one independent variable, CHP-1d, with the addition of the sampling points for its heat generation - $\hat{\text{H}}_{i,j}$.\\
Table \ref{Tab: CHP2d Data} gives the technical parameters for the CHP unit with two degrees-of-freedom, CHP-2d. Given that these unit has two independent variables, fuel consumption - $f_{i,t}$, and valve opening - $o_{i,t}$, the sampling points for its electric and thermal energy generation must be indexed on both sampling sets $J$ and $K$. Therefore, in Table \ref{Tab: CHP2d Data}, the sampling points for the electric and thermal generation, $\hat{P}_{i,j,k}$ and $\hat{H}_{i,j,k}$, are provided based on the sampling points for the fuel and valve opening, $\hat{F}_{i,j}$ and $\hat{O}_{i,k}$.
\begin{table}[h]
\centering
\caption{Technical parameters of the electric units}
\label{Tab: EP Data}
\begin{tabular}{lrrrrr}
\toprule
\toprule
Name  & Gen 1   & Gen 2   & Gen 3  & Gen 4 \\
\midrule
Bus   & 1       & 8       & 11     & 13    \\
UR    & 288     & 144     & 54     & 72    \\
DR    & 275     & 110     & 45     & 80    \\
$\hat{F}_{i,1}$ & 888.9   & 500.0   & 200.0  & 275.9 \\
$\hat{F}_{i,2}$ & 1 555.6 & 875.0   & 350.0  & 482.8 \\
$\hat{F}_{i,3}$ & 2 222.2 & 1 250.0 & 500.0  & 689.7 \\
$\hat{P}_{i,1}$ & 288     & 144     & 54     & 72    \\
$\hat{P}_{i,2}$ & 543.2   & 271.6   & 101.85 & 135.8 \\
$\hat{P}_{i,3}$ & 800     & 400     & 150    & 200   \\
\bottomrule
\bottomrule
\end{tabular}
\end{table}
\begin{table}[]
\centering
\caption{Technical parameters of the CHP unit with 1 independent variable}
\label{Tab: CHP1d Data}
\resizebox{\columnwidth}{!}{%
\begin{tabular}{lcccccccccccccc}
\toprule
\toprule
\multicolumn{12}{c}{CHP-1d} \\
\midrule
Bus & $UR$ & $DR$ & $\hat{F}_1$ & $\hat{F}_2$ & $\hat{F}_3$ & $\hat{P}_1$ & $\hat{P}_2$ & $\hat{P}_3$ & $\hat{H}_1$ & $\hat{H}_2$ & $\hat{H}_3$ \\
\midrule
5   & 50 & 80 & 408.39      & 687.67      & 982.79      & 60.20       & 162.19      & 300.00      & 282.94      & 447.39      & 572.73     \\
\bottomrule
\bottomrule
\end{tabular}
}
\end{table}
\begin{table}[]
\caption{Technical parameters of the CHP unit with 2 independent variables}
\label{Tab: CHP2d Data}
\resizebox{\columnwidth}{!}{%
\begin{tabular}{cccccccc|ccc}
\toprule
\toprule
\multicolumn{11}{c}{CHP-2d} \\
\midrule
Bus                & $UR$                  & $DR$                  & $\hat{F}_j$ & $\hat{O}_k$ & \multicolumn{3}{c}{$\hat{P}_{j,k}$} & \multicolumn{3}{c}{$\hat{H}_{j,k}$} \\
\midrule
\multirow{3}{*}{5} & \multirow{3}{*}{50} & \multirow{3}{*}{50} & 457.8       &\multicolumn{1}{c|}{0.0}         & 230.1      & 190.4      & 150.3     & 0.0       & 79.6       & 159.7      \\
                   &                     &                     & 658.2       &\multicolumn{1}{c|}{40}           & 362.2      & 323.0      & 283.2     & 0.0       & 102.1      & 204.7      \\
                   &                     &                     & 871.6       &\multicolumn{1}{c|}{80}           & 500.0      & 462.7      & 425.0     & 0.0       & 124.1      & 248.4   \\
\bottomrule
\bottomrule
\end{tabular}
}
\end{table}
\subsection{Thermal system}
\begin{figure}[htb]
\centering
\includegraphics[width=0.85\linewidth]{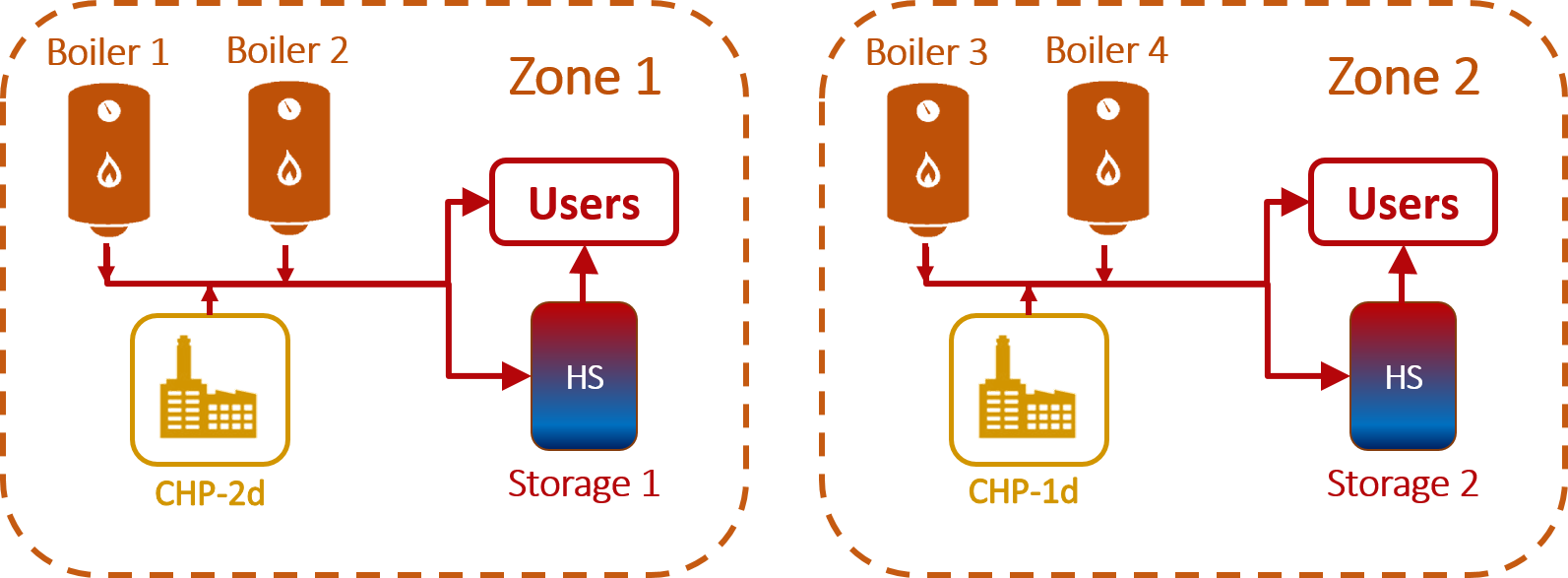}
\label{Fig: i^max_cha}
\caption{Configuration of thermal zones}
\label{fig1}
\end{figure}

The heating system is set up by assigning two thermal zones to the electrical system with a given heat load for each period. Each zone accounts two boilers, one CHP and a zonal thermal energy storage. The thermal capacity of the CHPs is taken to be 45~\% of the maximum load and the boilers act as a back-up to fulfill the entire load according to the operational practice of district heating networks \cite{Baccino2014a}. The topology of the heating system is shown in Fig. \ref{fig1}. Table \ref{Tab: HT Data} presents the technical parameters for the thermal generation units. Their serving zone, their sampling points for fuel consumption - $\hat{F}_{i,j}$ and thermal energy generation - $\hat{H}_{i,j}$ are presented.\\

\begin{table}[h]
\centering
\caption{Technical parameters of the thermal units}
\label{Tab: HT Data}
\begin{tabular}{lrrrr}
\toprule
\toprule
Name        & Boiler 1 & Boiler 2 & Boiler 3 & Boiler 4 \\
\midrule
Zone        & 1        & 1        & 2        & 2        \\
$\hat{F}_{i,1}$ & 27.2     & 22.7     & 65.2     & 51.1     \\
$\hat{F}_{i,2}$ & 149.5    & 125.0    & 358.7    & 281.3    \\
$\hat{F}_{i,3}$ & 271.7    & 227.3    & 652.2    & 511.4    \\
$\hat{H}_{i,1}$ & 25.0     & 20.0     & 60.0     & 45.0     \\
$\hat{H}_{i,2}$ & 137.5    & 110.0    & 330.0    & 247.5    \\
$\hat{H}_{i,3}$ & 250.0    & 200.0    & 600.0    & 450.0   \\
\bottomrule
\bottomrule
\end{tabular}
\end{table}
\subsection{Energy demand}
\begin{figure}[htb]
\centering
\resizebox{0.85\columnwidth}{!}{\input{energy_load.tex}}
\caption{Total electricity demand and zonal thermal energy request for a typical day}
\label{fig2}
\end{figure}
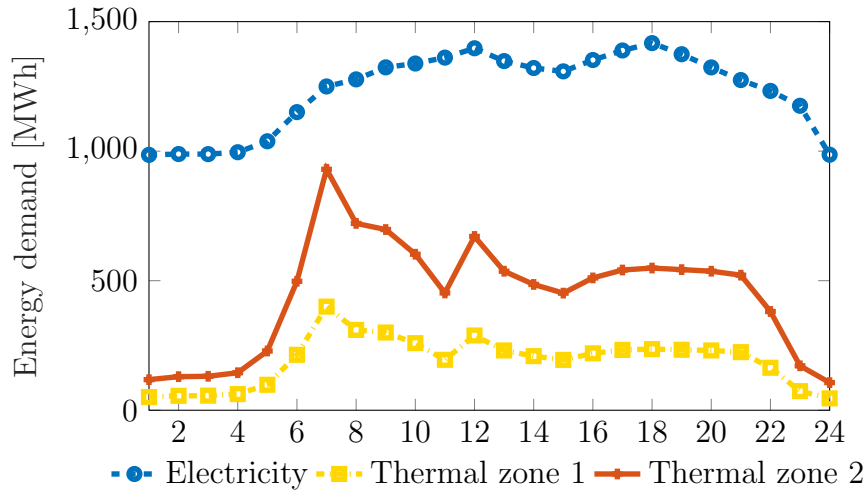

The load profile for the electric network and the thermal zones for a typical day are depicted in Fig. \ref{fig2}. The electric load per node has been normalized based on the energy demand for a typical winter day while keeping the nodal proportions for energy request present in the original IEEE 30-bus test case. The same approach has been implemented for designation of the thermal loads \cite{Baccino2014a}. The load data can be found in the online dataset \cite{Gonzalez-CastellanosData2018}.\\

\subsection{Implementation notes}
Given the MILP formulation of the model, it is possible to solve it by employing an  \enquote{off-the-shelf solver} to which the optimization data, variables and constraints are passed. The simulations are performed using the modeling software GAMS 24.9 \cite{GamsSoftware2013} with Gurobi 7.5 \cite{gurobi} as a MILP solver, on a computer with an Intel core i5-7200 @ 2.5GHz and 8 GB RAM.\\

\subsection{Results and discussion}
\subsubsection{Case 1: base case}
The constraints described in Section \ref{S: MathModel} are used to model the joint energy system. This is done to set a reference for the possible topologies and methodologies employed for the analysis of the combined energy system.\\
Figure \ref{fig4} shows the scheduling of the units which is characterized by the intensive use of the CHP units for both electrical and thermal dispatch in which they provide 59~\% and 96~\% of the requested electric and thermal energies respectively. In both thermal zones the economic benefits of employing CHP units are evidenced. Even though the boilers have a high first principle thermodynamics efficiency, the comparison must be performed by assessing the amount of fuel needed to generate heat (boiler) and power (electric generator) separately, against the CHP solution.\\
The thermal load is satisfied by the CHPs depending on their capacities and stored energies. The smart usage of thermal storage and CHPs can be seen explicitly at the sixth hour when the peak demand occurs. At this hour even though the demand highly exceeds the thermal generation limits of the units CHP-2d (orange bar) and CHP-1d (red bar) (Figures \ref{fig4b} and \ref{fig4c}) the use of the stored energy (green bar) removes the need of using the auxiliary boilers.\\
\begin{figure}
\centering\resizebox{0.65\columnwidth}{!}{
\subfloat[Electricity]{\input{ED1.tex}}%
\label{fig4a}}
\hfil
\resizebox{0.65\columnwidth}{!}{
\subfloat[Thermal zone 1\label{fig4b}]{\input{zone1_1.tex}}}%
\hfil
\resizebox{0.65\columnwidth}{!}{
\subfloat[Thermal zone 2\label{fig4c}]{\input{zone2_1.tex}}}%
\caption{Optimal schedules for case 1.}
\label{fig4}
\end{figure}
When the heat is neither consumed nor stored or lost, it is dissipated to the environment. The dissipated heat is only zero for zone 1. Analyzing the first four hours of the thermal dispatch at the second zone, Figure \ref{fig4b}, the thermal load (dashed line) is significantly lower than the heat production from the CHP-1d (red bar), which leads to heat dissipation during those hours (solid orange line). The economic benefits of dissipating excess heat in the second zone, instead of using the boilers, is a consequence of the higher combined efficiency of the CHP units when compared to the use of electric generators with boilers.\\
Table \ref{tab4} summarizes the main MILP model features and results obtained for this and the following test cases.
\begin{table}
  \centering
  \caption{Main model features and optimization results}
    \begin{tabular}{p{13em}p{3em}p{3em}p{3em}p{3em}}
    \toprule
    \toprule
    \multicolumn{1}{c}{} & Case 1 & Case 2 & Case 3 & Case 4 \\
    \midrule
    Number of binary variables & 1 608 & 1 608 & 1 608 & 1 608 \\
    Total number of variables & 6 126 & 6 078 & 6 126 & 6 126 \\
    Relative MILP gap [\%] & 5E-03 & 2E-03 & 9E-03 & 0 \\
    Computational time [s] & 10.14 & 2.13 & 1.79 & 2.67 \\
    Total cost [\euro]  & 5 585 & 5 660 & 5 403 & 6 201 \\
    Cost variation [\%] & - & 1.34 & -3.27 & 11.02 \\
    $\text{CO}_2$ emissions [ton]  & 7 083 & 7 188 & 6 974 & 7 674 \\
    \bottomrule
    \bottomrule
    \end{tabular}%
  \label{tab4}%
\end{table}%

\subsubsection{Case 2: operation without energy storage}
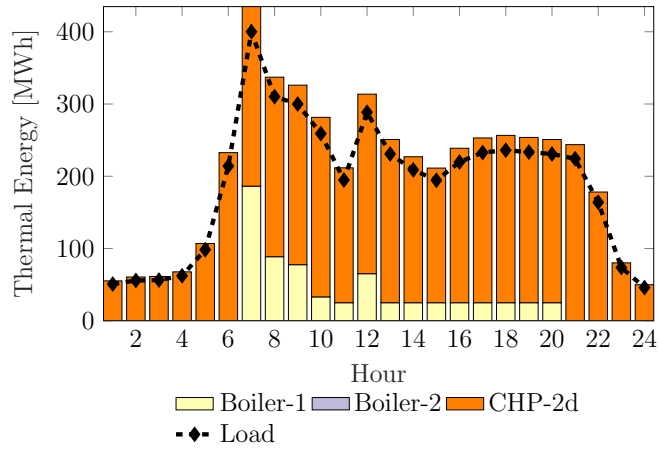
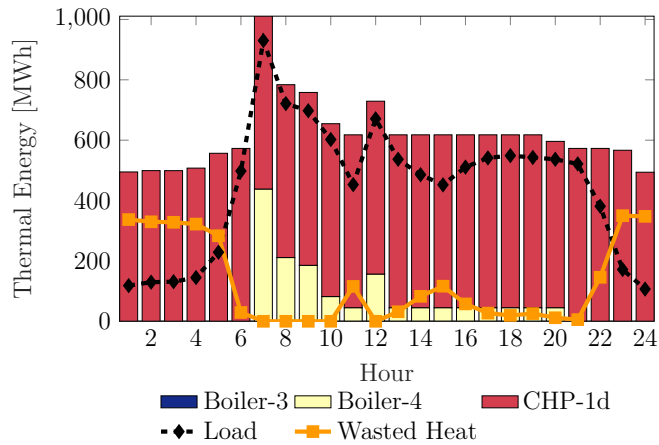
\begin{figure}[h]
\centering\resizebox{0.65\columnwidth}{!}{
\subfloat[Thermal zone 1\label{fig5a}]{\input{zone1_2.tex}}%
}
\hfil
\resizebox{0.65\columnwidth}{!}{
\subfloat[Thermal zone 2\label{fig5b}]{\input{zone2_2.tex}}%
}
\caption{Optimal schedules for case 2.}
\label{fig5}
\end{figure}

Unlike case 1, this case does not include the use of thermal energy storage. The rest of the system is modeled as in the base case.\\
As seen in Figure \ref{fig5}, the effect of a lack of storage is different at each zone. In both thermal zones, the use of boilers increases to aid the CHP units when the demand exceeds their rated capacity, which is evident in the hours following the demand peak at the sixth hour. For this case the boilers are used for 14 hours after the peak, while in the first one (Figure \ref{fig4}) they are only used for the three and four hours that follow the thermal peak at zones 1 and 2, respectively.
In the first zone of the first case it is possible to entirely satisfy the demand during the first 2 hours from the stored energy, whereas in the thermal zone 2 the heat is being dissipated during the same hours, as discussed above. Moreover, when the energy storage is used the amount of dissipated heat is reduced by 29~\% in the second zone. \\
The more flexible operation of the CHP with 2 degrees-of-freedom and a smart use of the energy storage system allows the reduction of wasted heat, thereby increasing the plant efficiency and primary energy savings. This type of operation is desired under the European Union policy framework, since it would allow the plants to qualify as high efficiency CHP systems and deserve obtaining financial incentives \cite{Bischi2017a}.\\
As expected, the use of energy storage reduces the operational costs of the system by about 1.3~\% as compared with the first case, due to the additional cost of operating the boilers to fulfill the demand peaks.\\

\subsubsection{Case 3: unconstrained network operation}
\begin{figure*}
\centering
\def\svgwidth{0.95\columnwidth}
\subfloat[Case 1\label{fig6a}]{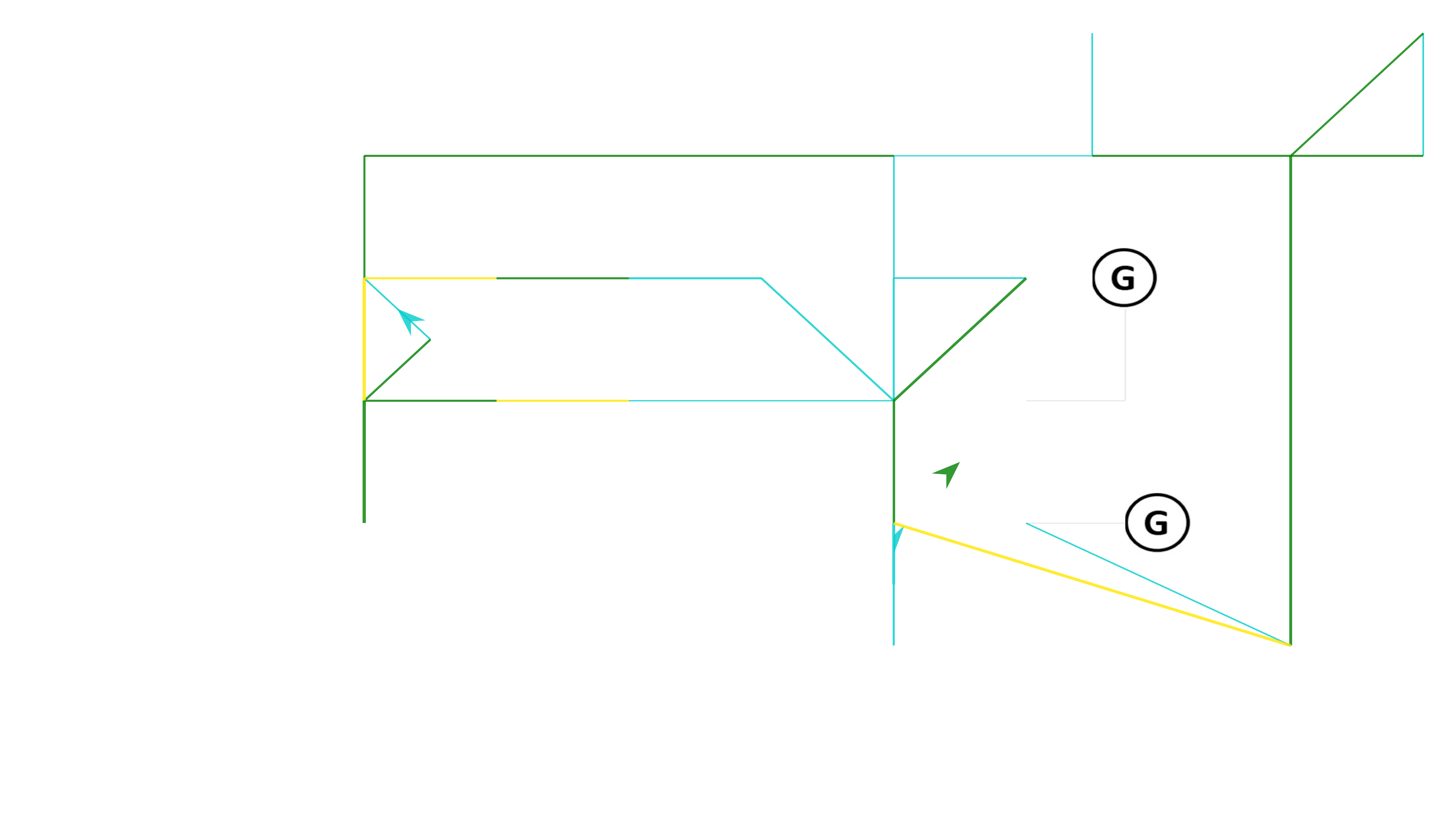}%
\hfil
\def\svgwidth{0.95\columnwidth}
\subfloat[Case 3\label{fig6b}]{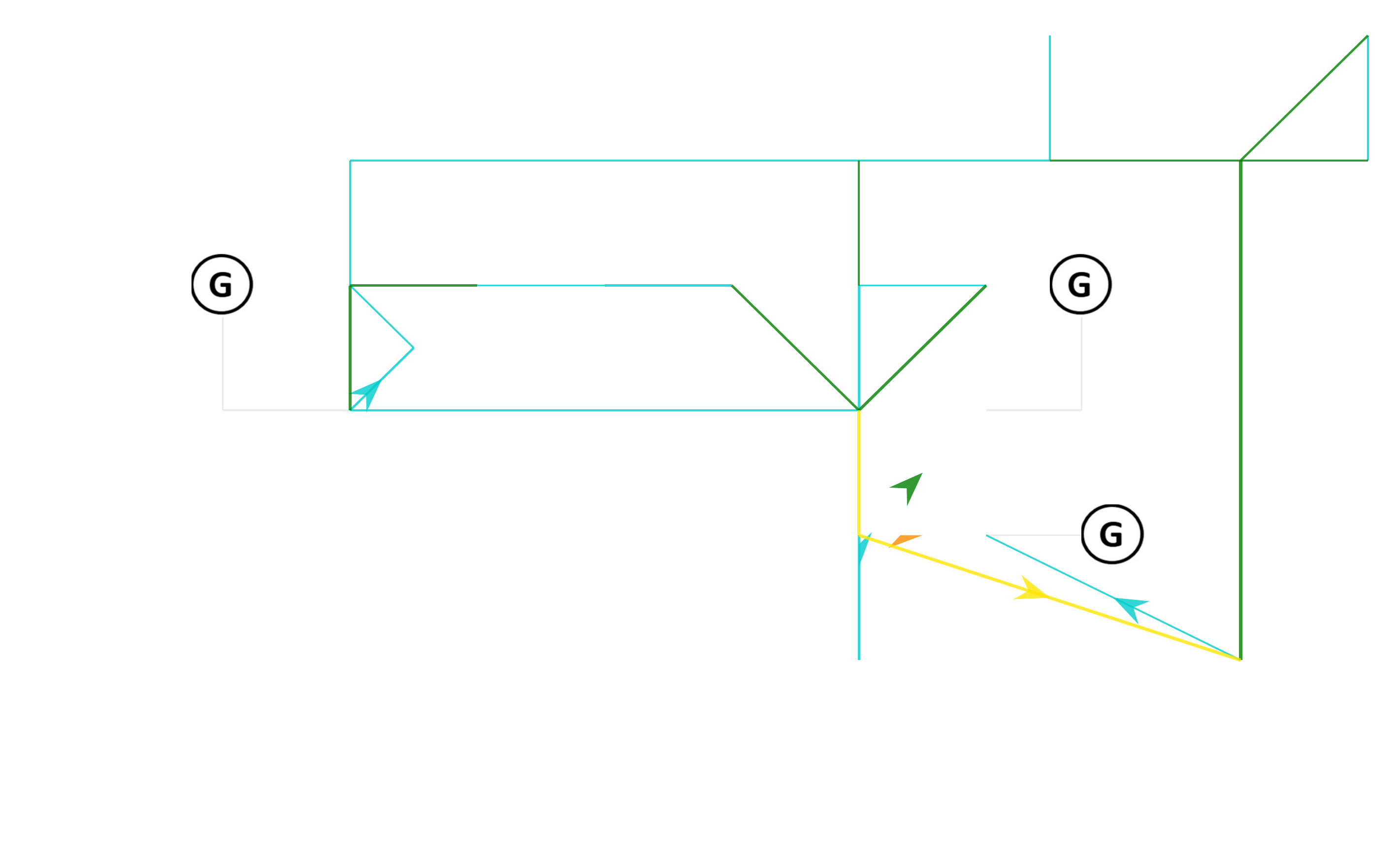}%
\caption{Power flow in the electric network on the most loaded hour. In red, overloaded lines, with a power flow higher than the STR. Orange, yellow, green and blue lines are loaded between 75-100~\%, 50-75~\%, 25-50~\% and below 25~\% of their STR, respectively. The gray lines indicate a lack of power input from the units.}
\label{fig6}
\end{figure*}

This case shows the impact of the network constraints on the unit commitment, especially on the CHP units that are closer to the overloaded lines. For this, an analysis of hour 18, during which the demand is the highest, is presented. The generation units as well as the thermal energy storage capacity are same as that of the base case.\\
The removal of the STR constraint \eqref{Eq: STR} makes the system operating at a price which is 3.3~\% lower than that of the first case. This operation, with the removal of the STR constraint, leads to the overloading of the line between nodes 2 and 5; marked in Figure \ref{fig6a} with a red color.\\
In the network-constrained operation, case 1, the generator 4 is committed for the hours surrounding the load peaks, 11 to 20. The introduction of this generator allows the decongestion of line 2-5, by serving more locally the demand present at nodes 12-22. When compared to case 3, without STR, the case 1 reduces the energy generated by the CHP-1d and the generator 1 by 18 and 16~\%, respectively.\\
Even though the CHP-2d is directly connected to the line 2-5, its total energy generation, thermal and electric, does not change. A more general sense would be to reduce the electric generation in one or both of the CHP units, connected at both ends of line 2-5, in order to relieve the overloading and to turn on a-priory more expensive generator. The only change in the operation of the CHP-2d, is presented in a change in the use of the storage system in the zone 1, to change the amount of heat and electricity produced around peak hours. By manipulating the \enquote{degrees-of-freedom} of the CHP-2d, it is possible to control the power flow in the electric system via the modulation of its energy fed-in; as a result, improving the operational flexibility of the system against congestion \cite{Ulbig2017}.\\

\subsubsection{Case 4: decoupled electric and thermal operation}
A decoupled operation is assumed, in which the thermal and electric energy systems are operated by independent entities. The thermal energy system has a dispatch priority, thereby prioritizing comfort of the users. For the thermal unit scheduling, the zones operating with CHPs having 2 degrees-of-freedom would make use of them as few as possible, favoring the use of boilers in order to sell more electricity in the market. Those with CHP units having one degree-of-freedom would pursue the opposite approach, otherwise, these units cannot produce electricity to be sold. With the results from the thermal dispatch i.e. electric generation of the CHP units, the electric energy dispatch is calculated for the fulfillment of the remaining load.\\
This case analyzes the impact of the application of a disjoint approach for the economic dispatch of the energy system. \\
As seen in Table \ref{tab4}, the implementation of a decoupled economic dispatch results in an increase of the operational costs by 11~\%. The cost increment is due to:
\begin{itemize}
\item non-inclusion of the CHP-2d in the thermal dispatch of the first heating zone. Even though its combined efficiency is reduced when it is producing both heat and electricity, it is still higher than the alternative combination of boilers and generators, as seen in the scheduling of the first case in Figure \ref{fig4}. The scheduling of this unit in the first case is a clear situation where the economic optimum differs from the thermodynamic one;
\item the reduction of the energy output of CHP-1d by 15~\% as compared to the first case. This is a consequence of operation of the heating zone 2 where no heat is dissipated. This contributes to an increase in the operation cost of the joint system as seen in the first case, since it is economically favorable to dissipate heat in that zone in order to reduce the costs associated with the electricity generation from more expensive units.
\end{itemize}
The joint operation of the thermal and electric system, case 1 \textit{vs.} case 3, allows to reap economic benefits, as well as increase the system flexibility and operational efficiency; by reducing the dissipated heat.\\

\section{Conclusion} \label{S: Conclusion}
A detailed MILP model has been developed to address a mixed integer non-linear problem for the combined heat and power (CHP) unit commitment with multiple heating zones, thermal energy storage and electric network constraints. Four case studies addressing the impact of the co-generators flexibility, i.e. one or two degrees-of-freedom, the presence of thermal energy storage, transmission constraints and the joint operation of the thermal and electric systems are considered.\\
An analysis of the evaluated test cases allows to come to the following conclusions: 
\begin{itemize}
    \item A flexible operation of CHPs with two-degrees of freedom results in district heating zones without heat dissipation. Such units allow a smarter energy management by enabling for some hours the fulfillment of the thermal load entirely from the thermal storage system, freeing the generation units to produce more electricity depending on the electric load and price profile. This operation improves the revenue possibilities of a CHP plant.
    \item The use of thermal energy storage increases the primary energy savings of the CHP system. The amount of dissipated heat is decreased by a quarter when the storage is employed. Additionally, the storage system reduces the need of using auxiliary boilers by more than 60~\% for the considered case. Resulting in economic benefits related to the savings in the operation and maintenance costs of these units.
    \item By considering a joint electricity and heat system when scheduling network-constrained CHP systems, the economic benefits of employing co-generation units have been confirmed by a power flow redistribution that allows them to operate at high generation rates, in spite of being directly connected to the overloaded lines. Such flow redistribution increased the system operational cost by about 3~\% for the considered case, which is considerably low in comparison to the value of lost load in northern countries, i.e. the amount that the end users would be willing to pay to guarantee an uninterrupted energy supply.
\end{itemize}

In systems with variable renewable power generation there is a need to compensate for their sudden fluctuations. This can be done with the use of flexible CHP-plants and heat storage, which would allow to balance the electricity and heat production, thereby minimizing the generation costs. Such a balance could be obtained with the  storing of excess heat, in comparison to the thermal demand, produced by the CHP units while the renewable resources are producing electricity. Once the renewable generation varies, the CHP units can adjust their generation output to serve the mismatched electric load.\\
\section*{Acknowledgement}
This work was supported by Skoltech NGP Program (Skoltech-MIT joint project).

\section*{References}
\bibliographystyle{model1-num-names}
\bibliography{mendeley_v2.bib}

\end{document}

%% file: Pack.tex
\usepackage{longtable} 
\usepackage{csquotes} 
\usepackage{eurosym} 
\usepackage{IEEEtrantools} 
\usepackage{mathtools} 
\usepackage{hyperref}
\usepackage{tikz}
\usetikzlibrary{plotmarks}
\usepackage{pgfplots}
\usepackage{booktabs} 
\usepackage{subfig}
\usepackage{graphicx}
\usepackage{caption} 
\usepackage{tablefootnote} 
\usepackage{multirow} 

%% file: energy_load.tex
%
%
\definecolor{mycolor1}{rgb}{1.00000,0.84314,0.00000}%
\definecolor{mycolor2}{rgb}{0.85000,0.32500,0.09800}%
\definecolor{mycolor3}{rgb}{0.00000,0.44706,0.74118}%
\begin{tikzpicture}

\begin{axis}[%
width=3.5in,
height=2in,
at={(1.339in,0.993in)},
scale only axis,
xmin=1,
xmax=24,
xlabel style={font=\color{white!15!black}},
xlabel={Hour},
ymin=0,
ymax=1500,
ylabel style={font=\color{white!15!black}},
ylabel={Energy demand [MWh]},
axis background/.style={fill=white},
legend style={draw=none,legend style={at={(0.5,-0.1)},anchor=north}},
legend columns=-1
]

\addplot [color=mycolor3, dashed, line width=2.0pt, mark=o, mark options={solid, mycolor3}]
  table[row sep=crcr]{%
1	985.637793455949\\
2	988.893139214007\\
3	988.527153596237\\
4	995.885390753504\\
5	1038.12783600451\\
6	1150.69730707013\\
7	1249.72531027826\\
8	1277.03939480445\\
9	1323.44251865748\\
10	1338.27456737762\\
11	1361.02346288216\\
12	1397.15972703669\\
13	1347.69388143496\\
14	1320.91914413496\\
15	1307.95554830553\\
16	1351.7582480323\\
17	1388.60722102144\\
18	1417\\
19	1374.12189551813\\
20	1323.23063224719\\
21	1274.43897067819\\
22	1232.50472384215\\
23	1175.12203145603\\
24	986.292715087748\\
};

\addlegendentry{Electricity}\addplot [color=mycolor1, dashdotted, line width=2.0pt, mark=square, mark options={solid, mycolor1}]
  table[row sep=crcr]{%
1	50.82\\
2	55.74\\
3	56.39\\
4	62.3\\
5	98.36\\
6	214.1\\
7	400\\
8	310.16\\
9	300\\
10	259.02\\
11	194.75\\
12	288.52\\
13	230.82\\
14	208.85\\
15	194.43\\
16	219.67\\
17	232.79\\
18	236.07\\
19	233.44\\
20	230.82\\
21	224.26\\
22	163.93\\
23	73.77\\
24	45.9\\
};
\addlegendentry{Thermal zone 1}

\addplot [color=mycolor2, line width=2.0pt, mark=+, mark options={solid, mycolor2}]
  table[row sep=crcr]{%
1	118.16\\
2	129.59\\
3	131.11\\
4	144.84\\
5	228.69\\
6	497.78\\
7	930\\
8	721.13\\
9	697.5\\
10	602.21\\
11	452.8\\
12	670.82\\
13	536.66\\
14	485.58\\
15	452.04\\
16	510.74\\
17	541.23\\
18	548.85\\
19	542.75\\
20	536.66\\
21	521.41\\
22	381.15\\
23	171.52\\
24	106.72\\
};
\addlegendentry{Thermal zone 2}

\end{axis}
\end{tikzpicture}%

%% file: ED1.tex
%
%
\definecolor{mycolor1}{rgb}{0.00000,0.44706,0.74118}%
\definecolor{mycolor2}{rgb}{0.59608,0.30588,0.63922}%
\definecolor{mycolor3}{rgb}{0.69804,0.87451,0.54118}%
\definecolor{mycolor4}{rgb}{0.30196,0.68627,0.29020}%
\definecolor{mycolor5}{rgb}{0.89412,0.10196,0.10980}%
\definecolor{mycolor6}{rgb}{1.00000,0.49804,0.00000}%
\begin{tikzpicture}

\begin{axis}[%
width=3.5in,
height=2in,
at={(0.88in,1.939in)},
scale only axis,
point meta min=0,
point meta max=1,
xmin=1,
xmax=24,
xlabel style={font=\color{white!15!black}},
xlabel={Hour},
ymin=0,
ymax=613.244567377622,
ylabel style={font=\color{white!15!black}},
ylabel={Electric Energy [MWh]},
axis background/.style={fill=white},
legend style={at={(0.5,-0.2)}, anchor=north, legend columns=4, legend cell align=center, align=left, fill=none, draw=none}
]
\addplot [color=mycolor1, line width=2.0pt, mark=square*, mark options={solid, fill=mycolor1, mycolor1}]
  table[row sep=crcr]{%
1	288\\
2	288\\
3	288\\
4	288\\
5	288\\
6	422.220202125256\\
7	524.695310278265\\
8	552.009394804453\\
9	598.412518657477\\
10	613.244567377622\\
11	563.993462882155\\
12	551.966461800935\\
13	550.663881434961\\
14	517.397765680433\\
15	503.351585348021\\
16	554.728248032299\\
17	563.646710404053\\
18	524.870389590995\\
19	561.992285109123\\
20	543.2\\
21	547.989187335423\\
22	486.199695003457\\
23	429.422566212772\\
24	288\\
};
\addlegendentry{Gen. 1}

\addplot [color=mycolor2, line width=2.0pt, mark=triangle*, mark options={solid, rotate=180, fill=mycolor2, mycolor2}]
  table[row sep=crcr]{%
1	0\\
2	0\\
3	0\\
4	0\\
5	0\\
6	0\\
7	0\\
8	0\\
9	0\\
10	0\\
11	0\\
12	0\\
13	0\\
14	0\\
15	0\\
16	0\\
17	0\\
18	0\\
19	0\\
20	0\\
21	0\\
22	0\\
23	0\\
24	0\\
};
\addlegendentry{Gen. 2}

\addplot [color=mycolor3, line width=2.0pt, mark=triangle*, mark options={solid, fill=mycolor3, mycolor3}]
  table[row sep=crcr]{%
1	0\\
2	0\\
3	0\\
4	0\\
5	0\\
6	0\\
7	0\\
8	0\\
9	0\\
10	0\\
11	0\\
12	0\\
13	0\\
14	0\\
15	0\\
16	0\\
17	0\\
18	0\\
19	0\\
20	0\\
21	0\\
22	0\\
23	0\\
24	0\\
};
\addlegendentry{Gen. 3}

\addplot [color=mycolor4, dashdotted, line width=2.0pt, mark=*, mark options={solid, fill=mycolor4, mycolor4}]
  table[row sep=crcr]{%
1	0\\
2	0\\
3	0\\
4	0\\
5	0\\
6	0\\
7	0\\
8	0\\
9	0\\
10	0\\
11	72\\
12	120.163265235755\\
13	72\\
14	72\\
15	72\\
16	72\\
17	99.9305106173841\\
18	167.099610409005\\
19	87.0996104090047\\
20	72\\
21	0\\
22	0\\
23	0\\
24	0\\
};
\addlegendentry{Gen. 4}

\addplot [color=mycolor5, dotted, line width=2.0pt, mark=triangle*, mark options={solid, rotate=180, fill=mycolor5, mycolor5}]
  table[row sep=crcr]{%
1	197.64380780269\\
2	200.899153560747\\
3	218.850820948463\\
4	271.8833253601\\
5	300\\
6	300\\
7	300\\
8	300\\
9	300\\
10	300\\
11	300\\
12	300\\
13	300\\
14	300\\
15	300\\
16	300\\
17	300\\
18	300\\
19	300\\
20	283.000632247189\\
21	300\\
22	300\\
23	300\\
24	220\\
};
\addlegendentry{CHP-1d}

\addplot [color=mycolor6, dashed, line width=2.0pt, mark=x, mark options={solid, fill=mycolor6, mycolor6}]
  table[row sep=crcr]{%
1	499.99398565326\\
2	499.99398565326\\
3	481.676332647774\\
4	436.002065393404\\
5	450.127836004513\\
6	428.477104944874\\
7	425.03\\
8	425.03\\
9	425.03\\
10	425.03\\
11	425.03\\
12	425.03\\
13	425.03\\
14	431.521378454526\\
15	432.603962957512\\
16	425.03\\
17	425.03\\
18	425.03\\
19	425.03\\
20	425.03\\
21	426.449783342771\\
22	446.305028838692\\
23	445.699465243259\\
24	478.292715087748\\
};
\addlegendentry{CHP-2d}

\end{axis}
\end{tikzpicture}%

%% file: zone1_1.tex
%
%
\definecolor{mycolor1}{rgb}{1.00000,1.00000,0.70196}%
\definecolor{mycolor2}{rgb}{0.74510,0.72941,0.85490}%
\definecolor{mycolor3}{rgb}{1.00000,0.49804,0.00000}%
\definecolor{mycolor4}{rgb}{0.10588,0.30980,0.20784}%
\begin{tikzpicture}

\begin{axis}[%
width=3.5in,
height=2in,
at={(2.6in,2.885in)},
scale only axis,
bar width=0.8,
xmin=0.6,
xmax=24.4,
xlabel style={font=\color{white!15!black}},
xlabel={Hour},
ymin=-144.49507567304,
ymax=439.130434782609,
ylabel style={font=\color{white!15!black}},
ylabel={Thermal Energy [MWh]},
axis background/.style={fill=white},
legend style={at={(0.5,-0.2)}, anchor=north, legend columns=3, legend cell align=left, align=center, fill=none, draw=none}
]
\addplot[ybar stacked, fill=mycolor1, draw=black, area legend] table[row sep=crcr] {%
1	57.8304303416694\\
2	61.9210731603288\\
3	0.36829221479515\\
4	0\\
5	0\\
6	0\\
7	190.680434782609\\
8	88.8830340264651\\
9	77.6369565217392\\
10	66.2057723623785\\
11	0\\
12	66.6621234567901\\
13	2.49555555555557\\
14	0\\
15	0\\
16	0\\
17	5.05029789361383\\
18	8.50572620027439\\
19	5.46212345679015\\
20	2.49555555555557\\
21	0\\
22	0\\
23	0\\
24	0\\
};
\addplot[forget plot, color=white!15!black] table[row sep=crcr] {%
0.6	0\\
24.4	0\\
};
\addlegendentry{Boiler-1}

\addplot[ybar stacked, fill=mycolor2, draw=black, area legend] table[row sep=crcr] {%
1	0\\
2	0\\
3	0\\
4	0\\
5	0\\
6	0\\
7	0\\
8	0\\
9	0\\
10	0\\
11	0\\
12	0\\
13	0\\
14	0\\
15	0\\
16	0\\
17	0\\
18	0\\
19	0\\
20	0\\
21	0\\
22	0\\
23	0\\
24	0\\
};
\addplot[forget plot, color=white!15!black] table[row sep=crcr] {%
0.6	0\\
24.4	0\\
};
\addlegendentry{Boiler-2}

\addplot[ybar stacked, fill=mycolor3, draw=black, area legend] table[row sep=crcr] {%
1	0.02\\
2	0.02\\
3	60.9331923985699\\
4	212.212466977387\\
5	165.559165102461\\
6	237.065217391304\\
7	248.45\\
8	248.45\\
9	248.45\\
10	248.45\\
11	248.45\\
12	248.45\\
13	248.45\\
14	227.010869565217\\
15	223.435408585889\\
16	248.45\\
17	248.45\\
18	248.45\\
19	248.45\\
20	248.45\\
21	243.760869565217\\
22	178.184782608696\\
23	180.184782608696\\
24	72.1850131080977\\
};
\addplot[forget plot, color=white!15!black] table[row sep=crcr] {%
0.6	0\\
24.4	0\\
};
\addlegendentry{CHP-2d}

\addplot[ybar, bar width=0.8, fill=mycolor4, draw=black, area legend] table[row sep=crcr] {%
1	57.8304303416694\\
2	61.9210731603288\\
3	0.36829221479515\\
4	-144.49507567304\\
5	-55.5049243269605\\
6	0\\
7	190.680434782609\\
8	9.31956521739125\\
9	0\\
10	-33.112294101509\\
11	-36.0453849108367\\
12	66.6621234567901\\
13	2.49555555555557\\
14	0\\
15	-12.0984520641504\\
16	-9.41525104208358\\
17	5.05029789361383\\
18	8.50572620027439\\
19	5.46212345679015\\
20	2.49555555555557\\
21	0\\
22	0\\
23	-100\\
24	-20.1197957167933\\
};
\addplot[forget plot, color=white!15!black] table[row sep=crcr] {%
0.6	0\\
24.4	0\\
};
\addlegendentry{Storage 1}

\addplot [color=black, dashed, line width=2.0pt, mark=diamond*, mark options={solid, fill=white, black}]
  table[row sep=crcr]{%
1	50.82\\
2	55.74\\
3	56.39\\
4	62.3\\
5	98.36\\
6	214.1\\
7	400\\
8	310.16\\
9	300\\
10	259.02\\
11	194.75\\
12	288.52\\
13	230.82\\
14	208.85\\
15	194.43\\
16	219.67\\
17	232.79\\
18	236.07\\
19	233.44\\
20	230.82\\
21	224.26\\
22	163.93\\
23	73.77\\
24	45.9\\
};
\addlegendentry{Load}

\end{axis}
\end{tikzpicture}%

%% file: zone2_1.tex
%
%
\definecolor{mycolor1}{rgb}{0.07843,0.16863,0.54902}%
\definecolor{mycolor2}{rgb}{1.00000,1.00000,0.70196}%
\definecolor{mycolor3}{rgb}{0.83529,0.24314,0.30980}%
\definecolor{mycolor4}{rgb}{0.10588,0.30980,0.20784}%
\definecolor{mycolor5}{rgb}{1.00000,0.60000,0.00000}%
\begin{tikzpicture}

\begin{axis}[%
width=3.5in,
height=2in,
at={(0.849in,1.391in)},
scale only axis,
bar width=0.8,
separate axis lines,
every outer x axis line/.append style={black},
every x tick label/.append style={font=\color{black}},
every x tick/.append style={black},
xmin=0.6,
xmax=24.4,
xlabel={Hour},
every outer y axis line/.append style={black},
every y tick label/.append style={font=\color{black}},
every y tick/.append style={black},
ymin=-321.024,
ymax=1020.97826086957,
ytick={-400, -200,    0,  200,  400,  600,  800, 1000},
ylabel={Thermal Energy [MWh]},
legend style={at={(0.5,-0.2)}, anchor=north, legend columns=3, legend cell align=left, align=left, fill=none, draw=none}
]
\addplot[ybar stacked, fill=black, area legend] table[row sep=crcr] {%
1	0\\
2	0\\
3	0\\
4	0\\
5	0\\
6	0\\
7	448.248260869565\\
8	16.7517391304348\\
9	0\\
10	0\\
11	0\\
12	160.13892345679\\
13	10.8315555555555\\
14	0\\
15	0\\
16	0\\
17	17.4693014809161\\
18	25.3721427530701\\
19	18.1901396497425\\
20	26.6362842384146\\
21	0\\
22	5.98000000000012\\
23	0\\
24	232.5\\
};
\addplot[forget plot, color=black] table[row sep=crcr] {%
0.6	0\\
24.4	0\\
};
\addlegendentry{Boiler-3}


\addplot[ybar stacked, fill=mycolor2, draw=black, area legend] table[row sep=crcr] {%
1	0\\
2	0\\
3	0\\
4	0\\
5	0\\
6	0\\
7	0\\
8	194.71938563327\\
9	185.422173913043\\
10	81.8460869565217\\
11	90.414392055824\\
12	0\\
13	0\\
14	0\\
15	0\\
16	0\\
17	0\\
18	0\\
19	0\\
20	0\\
21	0\\
22	0\\
23	0\\
24	0\\
};
\addplot[forget plot, color=black] table[row sep=crcr] {%
0.6	0\\
24.4	0\\
};
\addlegendentry{Boiler-4}

\addplot[ybar stacked, fill=mycolor3, draw=black, area legend] table[row sep=crcr] {%
1	479.635702561419\\
2	482.596482166055\\
3	498.923758781513\\
4	547.157516150025\\
5	572.73\\
6	572.73\\
7	572.73\\
8	572.73\\
9	572.73\\
10	572.73\\
11	572.73\\
12	572.73\\
13	572.73\\
14	572.73\\
15	572.73\\
16	572.73\\
17	572.73\\
18	572.73\\
19	572.73\\
20	557.268852375464\\
21	572.73\\
22	572.73\\
23	572.73\\
24	499.968952180538\\
};
\addplot[forget plot, color=black] table[row sep=crcr] {%
0.6	0\\
24.4	0\\
};
\addlegendentry{CHP-1d}

\addplot[ybar , fill=mycolor4, draw=black, area legend] table[row sep=crcr] {%
1	0\\
2	0\\
3	0\\
4	-143.976\\
5	-321.024\\
6	0\\
7	448.248260869565\\
8	16.7517391304348\\
9	0\\
10	0\\
11	-170.970479012346\\
12	160.13892345679\\
13	10.8315555555555\\
14	0\\
15	-71.6475985248577\\
16	-16.0202695972857\\
17	17.4693014809161\\
18	25.3721427530701\\
19	18.1901396497425\\
20	26.6362842384146\\
21	-5.98000000000012\\
22	5.98000000000012\\
23	-232.5\\
24	232.5\\
};
\addplot[forget plot, color=black] table[row sep=crcr] {%
0.6	0\\
24.4	0\\
};
\addlegendentry{Storage 2}

\addplot [color=black, dashed, line width=2.0pt,  mark=diamond*, mark options={solid, fill=white, black}]
  table[row sep=crcr]{%
1	118.16\\
2	129.59\\
3	131.11\\
4	144.84\\
5	228.69\\
6	497.78\\
7	930\\
8	721.13\\
9	697.5\\
10	602.21\\
11	452.8\\
12	670.82\\
13	536.66\\
14	485.58\\
15	452.04\\
16	510.74\\
17	541.23\\
18	548.85\\
19	542.75\\
20	536.66\\
21	521.41\\
22	381.15\\
23	171.52\\
24	106.72\\
};
\addlegendentry{Load}

\addplot [color=mycolor5, line width=2.0pt, mark=square*, mark options={solid, fill=mycolor5, mycolor5}]
  table[row sep=crcr]{%
1	323.104846356505\\
2	314.39876359277\\
3	327.899858078992\\
4	226.086994858023\\
5	0\\
6	19.8316000000001\\
7	0\\
8	0\\
9	0\\
10	0\\
11	0\\
12	0\\
13	0\\
14	41.3316000000001\\
15	8.9558093571309\\
16	0\\
17	0\\
18	0\\
19	0\\
20	0\\
21	0\\
22	151.1436\\
23	141.4916\\
24	562.501436006096\\
};
\addlegendentry{Wasted Heat}

\end{axis}
\end{tikzpicture}%

%% file: zone1_2.tex
%
%
\definecolor{mycolor1}{rgb}{1.00000,1.00000,0.70196}%
\definecolor{mycolor2}{rgb}{0.74510,0.72941,0.85490}%
\definecolor{mycolor3}{rgb}{1.00000,0.49804,0.00000}%
\begin{tikzpicture}

\begin{axis}[%
width=3.5in,
height=2in,
at={(2.6in,2.884in)},
scale only axis,
bar width=0.8,
xmin=0.6,
xmax=24.4,
xlabel style={font=\color{white!15!black}},
xlabel={Hour},
ymin=0,
ymax=434.782608695652,
ylabel style={font=\color{white!15!black}},
ylabel={Thermal Energy [MWh]},
axis background/.style={fill=white},
legend style={at={(0.5,-0.2)}, anchor=north, legend columns=3, legend cell align=left, align=left, fill=none, draw=none}
]
\addplot[ybar stacked, fill=mycolor1, draw=black, area legend] table[row sep=crcr] {%
1	0\\
2	0\\
3	0\\
4	0\\
5	0\\
6	0\\
7	186.332608695652\\
8	88.6804347826088\\
9	77.6369565217392\\
10	33.0934782608696\\
11	25\\
12	65.1586956521739\\
13	25\\
14	25\\
15	25\\
16	25\\
17	25\\
18	25\\
19	25\\
20	25\\
21	0\\
22	0\\
23	0\\
24	0\\
};
\addplot[forget plot, color=white!15!black] table[row sep=crcr] {%
0.6	0\\
24.4	0\\
};
\addlegendentry{Boiler-1}

\addplot[ybar stacked, fill=mycolor2, draw=black, area legend] table[row sep=crcr] {%
1	0\\
2	0\\
3	0\\
4	0\\
5	0\\
6	0\\
7	0\\
8	0\\
9	0\\
10	0\\
11	0\\
12	0\\
13	0\\
14	0\\
15	0\\
16	0\\
17	0\\
18	0\\
19	0\\
20	0\\
21	0\\
22	0\\
23	0\\
24	0\\
};
\addplot[forget plot, color=white!15!black] table[row sep=crcr] {%
0.6	0\\
24.4	0\\
};
\addlegendentry{Boiler-2}

\addplot[ybar stacked, fill=mycolor3, draw=black, area legend] table[row sep=crcr] {%
1	55.2391304347826\\
2	60.5869565217391\\
3	61.2934782608696\\
4	67.7173913043478\\
5	106.913043478261\\
6	232.717391304348\\
7	248.45\\
8	248.45\\
9	248.45\\
10	248.45\\
11	186.684782608696\\
12	248.45\\
13	225.891304347826\\
14	202.010869565217\\
15	186.336956521739\\
16	213.771739130435\\
17	228.032608695652\\
18	231.597826086957\\
19	228.739130434783\\
20	225.891304347826\\
21	243.760869565217\\
22	178.184782608696\\
23	80.1847826086956\\
24	49.8913043478261\\
};
\addplot[forget plot, color=white!15!black] table[row sep=crcr] {%
0.6	0\\
24.4	0\\
};
\addlegendentry{CHP-2d}

\addplot [color=black, dashed, line width=2.0pt,  mark=diamond*, mark options={solid, fill=white, black}]
  table[row sep=crcr]{%
1	50.82\\
2	55.74\\
3	56.39\\
4	62.3\\
5	98.36\\
6	214.1\\
7	400\\
8	310.16\\
9	300\\
10	259.02\\
11	194.75\\
12	288.52\\
13	230.82\\
14	208.85\\
15	194.43\\
16	219.67\\
17	232.79\\
18	236.07\\
19	233.44\\
20	230.82\\
21	224.26\\
22	163.93\\
23	73.77\\
24	45.9\\
};
\addlegendentry{Load}

\end{axis}
\end{tikzpicture}%

%% file: zone2_2.tex
%
%
\definecolor{mycolor1}{rgb}{0.07843,0.16863,0.54902}%
\definecolor{mycolor2}{rgb}{1.00000,1.00000,0.70196}%
\definecolor{mycolor3}{rgb}{0.83529,0.24314,0.30980}%
\definecolor{mycolor4}{rgb}{1.00000,0.60000,0.00000}%
\begin{tikzpicture}

\begin{axis}[%
width=3.5in,
height=2in,
at={(2.6in,2.884in)},
scale only axis,
bar width=0.8,
xmin=0.6,
xmax=24.4,
xlabel style={font=\color{white!15!black}},
xlabel={Hour},
ymin=0,
ymax=1010.86956521739,
ylabel style={font=\color{white!15!black}},
ylabel={Thermal Energy [MWh]},
axis background/.style={fill=white},
legend style={at={(0.5,-0.2)}, anchor=north, legend columns=3, legend cell align=left, align=left, fill=none, draw=none}
]
\addplot[ybar stacked, fill=mycolor1, draw=black, area legend] table[row sep=crcr] {%
1	0\\
2	0\\
3	0\\
4	0\\
5	0\\
6	0\\
7	0\\
8	0\\
9	0\\
10	0\\
11	0\\
12	0\\
13	0\\
14	0\\
15	0\\
16	0\\
17	0\\
18	0\\
19	0\\
20	0\\
21	0\\
22	0\\
23	0\\
24	0\\
};
\addplot[forget plot, color=white!15!black] table[row sep=crcr] {%
0.6	0\\
24.4	0\\
};
\addlegendentry{Boiler-3}

\addplot[ybar stacked, fill=mycolor2, draw=black, area legend] table[row sep=crcr] {%
1	0\\
2	0\\
3	0\\
4	0\\
5	0\\
6	0\\
7	438.139565217391\\
8	211.106956521739\\
9	185.422173913043\\
10	81.8460869565217\\
11	45\\
12	156.422173913043\\
13	45\\
14	45\\
15	45\\
16	45\\
17	45\\
18	45\\
19	45\\
20	45\\
21	0\\
22	0\\
23	0\\
24	0\\
};
\addplot[forget plot, color=white!15!black] table[row sep=crcr] {%
0.6	0\\
24.4	0\\
};
\addlegendentry{Boiler-4}

\addplot[ybar stacked, fill=mycolor3, draw=black, area legend] table[row sep=crcr] {%
1	494.738485747238\\
2	499.161929781227\\
3	499.022299255084\\
4	507.471693378252\\
5	556.612010639642\\
6	572.73\\
7	572.73\\
8	572.73\\
9	572.73\\
10	572.73\\
11	572.73\\
12	572.73\\
13	572.73\\
14	572.73\\
15	572.73\\
16	572.73\\
17	572.73\\
18	572.73\\
19	572.73\\
20	551.056548288768\\
21	572.73\\
22	572.73\\
23	566.632528939462\\
24	493.87148112\\
};
\addplot[forget plot, color=white!15!black] table[row sep=crcr] {%
0.6	0\\
24.4	0\\
};
\addlegendentry{CHP-1d}

\addplot [color=black, dashed, line width=2.0pt,  mark=diamond*, mark options={solid, fill=white, black}]
  table[row sep=crcr]{%
1	118.16\\
2	129.59\\
3	131.11\\
4	144.84\\
5	228.69\\
6	497.78\\
7	930\\
8	721.13\\
9	697.5\\
10	602.21\\
11	452.8\\
12	670.82\\
13	536.66\\
14	485.58\\
15	452.04\\
16	510.74\\
17	541.23\\
18	548.85\\
19	542.75\\
20	536.66\\
21	521.41\\
22	381.15\\
23	171.52\\
24	106.72\\
};
\addlegendentry{Load}

\addplot [color=mycolor4, line width=2.0pt, mark=square*, mark options={solid, fill=mycolor4, mycolor4}]
  table[row sep=crcr]{%
1	336.999406887459\\
2	329.638975398729\\
3	327.990515314678\\
4	322.033957907992\\
5	283.39304978847\\
6	29.1316000000001\\
7	0\\
8	0\\
9	0\\
10	0\\
11	115.5116\\
12	0\\
13	31.6516\\
14	82.7316000000001\\
15	116.2716\\
16	57.5716\\
17	27.0816\\
18	19.4616\\
19	25.5616\\
20	11.7120244256662\\
21	5.50160000000011\\
22	145.7616\\
23	349.781926624305\\
24	347.6417626304\\
};
\addlegendentry{Wasted Heat}

\end{axis}
\end{tikzpicture}%

%% file: PF1.pdf_tex
\begingroup%
  \makeatletter%
  \providecommand\color[2][]{%
    \errmessage{(Inkscape) Color is used for the text in Inkscape, but the package 'color.sty' is not loaded}%
    \renewcommand\color[2][]{}%
  }%
  \providecommand\transparent[1]{%
    \errmessage{(Inkscape) Transparency is used (non-zero) for the text in Inkscape, but the package 'transparent.sty' is not loaded}%
    \renewcommand\transparent[1]{}%
  }%
  \providecommand\rotatebox[2]{#2}%
  \ifx\svgwidth\undefined%
    \setlength{\unitlength}{1016.93492344bp}%
    \ifx\svgscale\undefined%
      \relax%
    \else%
      \setlength{\unitlength}{\unitlength * \real{\svgscale}}%
    \fi%
  \else%
    \setlength{\unitlength}{\svgwidth}%
  \fi%
  \global\let\svgwidth\undefined%
  \global\let\svgscale\undefined%
  \makeatother%
  \begin{picture}(1,0.58059335)%
    \put(0,0){\includegraphics[width=\unitlength,page=1]{PF1.pdf}}%
    \put(0.69790889,0.27532352){\makebox(0,0)[lb]{\smash{ 9}}}%
    \put(0,0){\includegraphics[width=\unitlength,page=2]{PF1.pdf}}%
    \put(0.03679225,0.21181636){\makebox(0,0)[lb]{\smash{ 1}}}%
    \put(0.03679225,0.12714014){\makebox(0,0)[lb]{\smash{ 2}}}%
    \put(0.14888467,0.22875159){\makebox(0,0)[lb]{\smash{ 3}}}%
    \put(0.24038876,0.1906473){\makebox(0,0)[lb]{\smash{ 4}}}%
    \put(0.62699329,0.12714014){\makebox(0,0)[lb]{\smash{ 5}}}%
    \put(0.58581654,0.22875159){\makebox(0,0)[lb]{\smash{ 6}}}%
    \put(0.62699329,0.16947825){\makebox(0,0)[lb]{\smash{ 7}}}%
    \put(0.69790889,0.22875159){\makebox(0,0)[lb]{\smash{ 8}}}%
    \put(0.58581654,0.27532352){\makebox(0,0)[lb]{\smash{ 10}}}%
    \put(0.76653696,0.27532352){\makebox(0,0)[lb]{\smash{ 11}}}%
    \put(0.2609771,0.27532352){\makebox(0,0)[lb]{\smash{ 12}}}%
    \put(0.14888467,0.27532352){\makebox(0,0)[lb]{\smash{ 13}}}%
    \put(0.30672914,0.33883068){\makebox(0,0)[lb]{\smash{ 14}}}%
    \put(0.21980035,0.38116879){\makebox(0,0)[lb]{\smash{ 15}}}%
    \put(0.33189277,0.27532352){\makebox(0,0)[lb]{\smash{ 16}}}%
    \put(0.42339678,0.27532352){\makebox(0,0)[lb]{\smash{ 17}}}%
    \put(0.33189277,0.39810402){\makebox(0,0)[lb]{\smash{ 18}}}%
    \put(0.42339678,0.39810402){\makebox(0,0)[lb]{\smash{ 19}}}%
    \put(0.51490087,0.39810402){\makebox(0,0)[lb]{\smash{ 20}}}%
    \put(0.69790889,0.39810402){\makebox(0,0)[lb]{\smash{ 21}}}%
    \put(0.58581654,0.38116879){\makebox(0,0)[lb]{\smash{ 22}}}%
    \put(0.24038876,0.48278023){\makebox(0,0)[lb]{\smash{ 23}}}%
    \put(0.60640488,0.48278023){\makebox(0,0)[lb]{\smash{ 24}}}%
    \put(0.76424935,0.48278023){\makebox(0,0)[lb]{\smash{ 25}}}%
    \put(0.74366093,0.56745644){\makebox(0,0)[lb]{\smash{ 26}}}%
    \put(0.86032858,0.48278023){\makebox(0,0)[lb]{\smash{ 27}}}%
    \put(0.88091706,0.10597109){\makebox(0,0)[lb]{\smash{ 28}}}%
    \put(0.97242115,0.44467595){\makebox(0,0)[lb]{\smash{ 29}}}%
    \put(0.97242115,0.56745644){\makebox(0,0)[lb]{\smash{ 30}}}%
    \put(0.03716927,0.34968809){\makebox(0,0)[lb]{\smash{\textbf{Gen. 1}}}}%
    \put(0.76934178,0.24456668){\makebox(0,0)[lb]{\smash{\textbf{Gen. 2}}}}%
    \put(0.74628802,0.41282443){\makebox(0,0)[lb]{\smash{\textbf{Gen. 3}}}}%
    \put(0.12944272,0.41232122){\makebox(0,0)[lb]{\smash{\textbf{Gen. 4}}}}%
    \put(0.57898793,0.0003349){\makebox(0,0)[lb]{\smash{\textbf{CHP-1d}}}}%
    \put(0.02955927,0.00060292){\makebox(0,0)[lb]{\smash{\textbf{CHP-2d}}}}%
  \end{picture}%
\endgroup%

%% file: PF3.pdf_tex
\begingroup%
  \makeatletter%
  \providecommand\color[2][]{%
    \errmessage{(Inkscape) Color is used for the text in Inkscape, but the package 'color.sty' is not loaded}%
    \renewcommand\color[2][]{}%
  }%
  \providecommand\transparent[1]{%
    \errmessage{(Inkscape) Transparency is used (non-zero) for the text in Inkscape, but the package 'transparent.sty' is not loaded}%
    \renewcommand\transparent[1]{}%
  }%
  \providecommand\rotatebox[2]{#2}%
  \ifx\svgwidth\undefined%
    \setlength{\unitlength}{818.70092344bp}%
    \ifx\svgscale\undefined%
      \relax%
    \else%
      \setlength{\unitlength}{\unitlength * \real{\svgscale}}%
    \fi%
  \else%
    \setlength{\unitlength}{\svgwidth}%
  \fi%
  \global\let\svgwidth\undefined%
  \global\let\svgscale\undefined%
  \makeatother%
  \begin{picture}(1,0.61255052)%
    \put(0,0){\includegraphics[width=\unitlength,page=1]{PF3.pdf}}%
    \put(0.69311639,0.28847961){\makebox(0,0)[lb]{\smash{ 9}}}%
    \put(0,0){\includegraphics[width=\unitlength,page=2]{PF3.pdf}}%
    \put(0.03653963,0.2215767){\makebox(0,0)[lb]{\smash{ 1}}}%
    \put(0.03653963,0.13237288){\makebox(0,0)[lb]{\smash{ 2}}}%
    \put(0.14786233,0.23941747){\makebox(0,0)[lb]{\smash{ 3}}}%
    \put(0.238738,0.19927579){\makebox(0,0)[lb]{\smash{ 4}}}%
    \put(0.62268774,0.13237288){\makebox(0,0)[lb]{\smash{ 5}}}%
    \put(0.58179368,0.23941747){\makebox(0,0)[lb]{\smash{ 6}}}%
    \put(0.62268774,0.17697479){\makebox(0,0)[lb]{\smash{ 7}}}%
    \put(0.69311639,0.23941747){\makebox(0,0)[lb]{\smash{ 8}}}%
    \put(0.58179368,0.28847961){\makebox(0,0)[lb]{\smash{ 10}}}%
    \put(0.76127314,0.28847961){\makebox(0,0)[lb]{\smash{ 11}}}%
    \put(0.25918503,0.28847961){\makebox(0,0)[lb]{\smash{ 12}}}%
    \put(0.14786233,0.28847961){\makebox(0,0)[lb]{\smash{ 13}}}%
    \put(0.30462287,0.35538243){\makebox(0,0)[lb]{\smash{ 14}}}%
    \put(0.21829098,0.39998434){\makebox(0,0)[lb]{\smash{ 15}}}%
    \put(0.32961368,0.28847961){\makebox(0,0)[lb]{\smash{ 16}}}%
    \put(0.42048936,0.28847961){\makebox(0,0)[lb]{\smash{ 17}}}%
    \put(0.32961368,0.41782511){\makebox(0,0)[lb]{\smash{ 18}}}%
    \put(0.42048936,0.41782511){\makebox(0,0)[lb]{\smash{ 19}}}%
    \put(0.51136503,0.41782511){\makebox(0,0)[lb]{\smash{ 20}}}%
    \put(0.69311639,0.41782511){\makebox(0,0)[lb]{\smash{ 21}}}%
    \put(0.58179368,0.39998434){\makebox(0,0)[lb]{\smash{ 22}}}%
    \put(0.238738,0.50702893){\makebox(0,0)[lb]{\smash{ 23}}}%
    \put(0.60224071,0.50702893){\makebox(0,0)[lb]{\smash{ 24}}}%
    \put(0.75900125,0.50702893){\makebox(0,0)[lb]{\smash{ 25}}}%
    \put(0.73855422,0.59623275){\makebox(0,0)[lb]{\smash{ 26}}}%
    \put(0.85442062,0.50702893){\makebox(0,0)[lb]{\smash{ 27}}}%
    \put(0.87486774,0.11007197){\makebox(0,0)[lb]{\smash{ 28}}}%
    \put(0.96574341,0.46688725){\makebox(0,0)[lb]{\smash{ 29}}}%
    \put(0.96574341,0.59623275){\makebox(0,0)[lb]{\smash{ 30}}}%
    \put(0.03424153,0.36661584){\makebox(0,0)[lb]{\smash{\textbf{Gen. 1}}}}%
    \put(0.76090863,0.25641873){\makebox(0,0)[lb]{\smash{\textbf{Gen. 2}}}}%
    \put(0.73852802,0.43259799){\makebox(0,0)[lb]{\smash{\textbf{Gen. 3}}}}%
    \put(0.12511721,0.43443951){\makebox(0,0)[lb]{\smash{\textbf{Gen. 4}}}}%
    \put(0.57067726,0.000416){\makebox(0,0)[lb]{\smash{\textbf{CHP-1d}}}}%
    \put(0.02576921,0.00042021){\makebox(0,0)[lb]{\smash{\textbf{CHP-2d}}}}%
  \end{picture}%
\endgroup%